\documentclass[12pt]{article}
\usepackage{amsmath, amsthm}
\usepackage{newpxmath, newpxtext}
\usepackage{graphicx}\graphicspath{{./}}
\usepackage{fullpage}
\usepackage{biblatex}
\usepackage{longtable}
\usepackage{lscape, listings}

\usepackage{hyperref}
\hypersetup{
    colorlinks=true,
    linkcolor=blue,
    filecolor=magenta,      
    urlcolor=red,
    }
\newtheorem{theorem}{Theorem}

\newcommand{\RR}{\mathbb{R}}
\newcommand{\ZZ}{\mathbb{Z}}
\newcommand{\lt}{<}
\newcommand{\gt}{>}
\newcommand{\dee}[1]{\mathrm{d}#1}
\newcommand{\abs}[1]{\left| #1 \right|}
\newcommand{\set}[1]{\left\{ #1 \right\} }

\newcommand{\pdiff}[2]{\frac{\partial #1}{\partial #2}}
\newcommand{\bess}[2]{J\left(#1,#2 \right)}
\newcommand{\mybess}[2]{\mathcal{J}\left(#1,#2 \right)}
\newcommand{\const}{\mathrm{const}}

\DeclareMathOperator{\claus}{C}

\title{The first 128 digits of an autoconvolution inequality}
\author{Andrew Rechnitzer\thanks{andrewr@math.ubc.ca}\\
Department of Mathematics\\
University of British Columbia
}
\begin{document}
\maketitle

\begin{abstract}
	Using rigorous high-precision floating point arithmetic we compute
	very tight rigorous bounds on the auto-convolution constant
	\[
		\nu_2^2 = \inf_f \|f \ast f\|_2^2 = \inf_f \int_{-1}^1 (f \ast f)^2
	\]
	where the infimum is taken over all unit mass functions \(f \in
	L^1(-1/2,1/2)\). This quantity arises in additive combinatorics, particularly
	in the study of Sidon sets.  Our bounds give the first 128 digits of \(\nu_2^2\), and so
	substantially improve previous bounds on this quantity due to White, Green
	and Martin \& O'Bryant.
\end{abstract}

\section{Introduction}
A subset \(A \subseteq \set{1,2,\dots, N}\) is a \(B_2[g]\)-set, when for any \(n \in \ZZ \) there are at most
\(g\) representations of the form \(n = a_1+a_2\) where \(a_1,a_2 \in A\), up to reordering the summands.
When \(g=1\) these are precisely the classical Sidon sets. A major problem in additive number theory is understanding
the behaviour of the quantity \(R_2[g](N)\), the size of the largest \(B_2[g]\)-set. This problem dates back to a question of
Sidon~\cite{sidon1932} and work of Erd\H{o}s and Tur{\'a}n~\cite{erdos1941}. The problem is readily generalised to \(h\) summands
giving \(B_h[g]\)-sets and \(R_h[g](N)\), and we refer the reader to \cite{obraynt2004} for a survey of Sidon sets and their generalisations.

A counting argument gives an upper bound\(R_h[g](N) \leq \sqrt{g h h! N}\), while the work of several authors~\cite{lindstrom2000, cilleruelo2000, johnston2022} gives a lower bound \(R_h[g](N) \geq (1-o(1))(gN)^{1/h}\). Estimating the limit
\begin{equation}
	\sigma_h(g) = \lim_{N \to \infty} \frac{R_h[g](N)}{(gN)^{1/h}}
\end{equation}
is an open problem. Indeed, the limit is only known to exist for classical Sidon sets where \(\sigma_2(1)= 1\). Cilleruelo, Ruzsa \& Vinuesa~\cite{cilleruelo2010} have shown that \(\displaystyle \lim_{g \to \infty} \sigma_2(g)\) exists and is expressed in terms of the quantity
\begin{equation}
	\nu_\infty = \inf_{f \in \mathcal{F}} \| f \ast f \|_\infty
\end{equation}
where the infimum is taken over all non-negative functions \(f \in
L^1(-1/2,1/2)\). This quantity was studied by Schinzel \&
Schmidt~\cite{schinzel2002} and linked to Sidon sets by several
authors~\cite{cilleruelo2008,martin2009}, but remains very stubbornly difficult
to compute. It is a statement of just how difficult the underlying problem is, that
the current\footnote{There has been a succession of results in
	recent years giving upper bounds by constructing step-functions via machine
	learning methods --- such as \cite{georgiev2025, yuksekgonul2026}; there may be even newer results.}
best lower bound~\cite{cloninger2017} and upper bound~\cite{yuksekgonul2026}
do not quite specify the second digit.

More generally, for any \(1 < p \leq \infty \) we then define
\begin{equation}
	\nu_p^p = \inf_{f \in \mathcal{F}} \| f\ast f\|_p^p =
	\int_{f \in \mathcal{F}} \left( \int_{-1}^1\left( \int_{-1/2}^{1/2} f(t) f(x-t) \dee{t} \right)^p \right).
\end{equation}
The question of computing this for any \(p\) appears as Problem~35
in~\cite{green100}, and remains unsolved for any \(p \geq 2\). When \(p=2\), the subject of this work,
the constant \(\nu_2\) connects, via work of Green~\cite{green2001} and White~\cite{white2022}, to bounds for \(\sigma_2(g)\) for small \(g\) via
\begin{equation}
	\sigma_2(g) \leq \frac{\sqrt{2-1/g}}{\nu_2}.
\end{equation}
Prior to this work, the strongest bounds on \(\nu_2^2\) are due to White~\cite{white2024}, improving on previous bounds
due to Green~\cite{green2001} and Martin \& O'Bryant~\cite{martin2007}:
\begin{align}
	0.574575 & \lt \nu_2^2 \lt 0.640733, \nonumber \\
	0.574636 & \lt \nu_2^2 \lt 0.574643.
\end{align}
The first set of bounds are due to Green (lower)~\cite{green2001} and Martin \& O'Bryant (upper)~\cite{martin2007},
while the second set of tighter bounds are due to White~\cite{white2024} --- these bounds give the first 4 digits.

In this work we compute bounds that give the first 128 digits --- see Theorem~\ref{thm main}. In Section~\ref{sec first ansatz}
we describe results of White~\cite{white2024} which converts the problem of minimising \(\|f\ast f\|_2\) into a quadratic
programming problem. By analysing the (near optimal) Fourier coefficients in White's work we were led to a first ansatz.
This allows us to translate the problem from one of minimising over a very large set of Fourier coefficients to one of
minimising over a much smaller set of parameters. This results in much tighter numerical bounds but we were unable to make them
completely rigorous. However, the resulting data led us to a second ansatz, described in Section~\ref{sec second ansatz} which we can then use, with some careful rigorous floating point summations, to obtain rigorous upper bounds.
These are then leveraged in Section~\ref{sec lower} to find rigorous lower bounds. Finally in Section~\ref{sec results} we give details of our rigorous floating point computations, state our main result Theorem~\ref{thm main}, and explore
our (near) optimal auto-convolution function.

\section{A first ansatz}\label{sec first ansatz}
As above, let \(\mathcal{F}\) denote the set of non-negative functions in \(L^1(-1/2,1/2)\), then
for \(f \in \mathcal{F}\) we define its Fourier transform as
\begin{equation}
	\hat{f}(k)  = \int_{-1/2}^{1/2} e^{-2\pi i k x} f(x) \dee{x} \qquad k \in \mathbb{Z}.
\end{equation}
White~\cite{white2024} shows that there is a unique optimal choice of \(f\) giving \(\nu_2\)
and that it be an even function, so that \(\hat{f}(-k)=\hat{f}(k)\).
Then extend \(f\) on \((-1/2,1/2)\) to a function \(F\) defined on \((-1,1)\)
and its transform \(\hat{F}\):
\begin{align}
	F(x)       & = \begin{cases}
		               f(x) & x \in (-1/2,1/2) \\
		               0    & \text{otherwise}
	               \end{cases}, \text{ and}                           \nonumber \\
	\hat{F}(k) & = \frac{1}{2} \int_{-1}^{1} e^{-\pi i k x} F(x) \dee{x}
	= \frac{1}{2} \int_{-1/2}^{1/2} e^{-\pi i k x} f(x) \dee{x}.
	\label{eqn_big_f_hat}
\end{align}
White (see \cite{white2024} Lemma~3.1) then establishes a relationship between \(\hat{f}\) and \(\hat{F}\) and rewrites \(\nu_2^2\)  as
\begin{align}
	\nu_2^2        & = \inf_f \mathcal{C}(f)                                     \qquad\qquad \text{where} \nonumber \\
	\label{eqn_c_by_big_f}
	\mathcal{C}(f) & = \| F\ast F\|_2^2  = 8 \sum_{k \in \ZZ} \abs{\hat{F}(k)}^4
	= \frac{1}{2} \sum_{k \in \ZZ} \abs{\hat{f}(k)}^4 + \frac{8}{\pi^4} \sum_{m \in \ZZ} \abs{ \sum_{k \in \ZZ} \frac{(-1)^k \hat{f}(k)}{2k-(2m+1)} }^4
\end{align}
By truncating the infinite sums in the above expression, White transforms the problem of minimising \(\mathcal{C}(f)\) into one of optimising
a finite set of coefficients \(\set{ \hat{f}(1), \hat{f}(2), \dots, \hat{f}(T) } \) for some large but finite \(T\), implicitly setting
\(\hat{f}(k)=0\) when \(|k| \gt T\). The resulting optimisation problem can be written as a quadratically constrained linear program, and
can be solved numerically using fairly standard software packages. That solution then gives rigorous upper and (with some extra work) lower bounds
for \(\nu_2^2\).

We obtained the near optimal values \(\hat{f}(k)\) from White~\cite{white_data}. They are alternating and decrease in magnitude fairly rapidly with \(k\).
Some simple fitting shows that these coefficients lie very nearly on a curve of the form
\begin{equation}
	(-1)^k \hat{f}(k) = \frac{a}{\sqrt{k}},
\end{equation}
where the constant \(a\) is approximately \(0.3\). Some more careful numerics suggests that
\begin{equation}
	\label{eqn_ansatz_1}
	(-1)^k\hat{f}(k) = \frac{1}{\sqrt{k}} \sum_{j=0}^{P-1} a_j k^{-j},
\end{equation}
for some modest value of \(P\), will give an even better fit.

We use the above as an ansatz to replace the problem of optimising directly over
the large set \(\set{\hat{f}(1),\dots, \hat{f}(T)}\), with one of optimising over a much smaller
set of coefficients \(\vec{a} = (a_0,a_1,\dots,a_P) \). With this ansatz we compute \(\mathcal{C}(f)\)
in two pieces. The first term
\begin{align}
	\sum_{k \in \ZZ} |\hat{f}(k)|^4 & = \frac{1}{2} + \sum_{i,j,k,\ell=0}^{P-1} a_i a_j a_k a_\ell \cdot \zeta(2+i+j+k+\ell),
\end{align}
where \(\zeta(s)\) can be computed to high precision, making this a simple polynomial in the \(a_j\).
Then
\begin{align}
	\frac{8}{\pi^4} \sum_{m \in \ZZ} \abs{ \sum_{k \in \ZZ} \frac{(-1)^k \hat{f
			}(k)}{2k-(2m+1)} }^4
	                      & = \frac{16}{\pi^4} \sum_{m \geq 1} \left(
	\frac{1}{2m-1} + 2 \sum_{k\geq 1} \frac{m \hat{f}(k) (-1)^k}{(2m-1)^2-4k^2}
	\right)^4                                                 \nonumber       \\
	                      & = \frac{16}{\pi^4} \sum_{m \in \ZZ}
	\left(
	\frac{1}{2m-1} + 2 \sum_p a_p \cdot E_{m,p}
	\right)^4                                                  \nonumber      \\
	\text{where }	E_{m,p} & = \sum_{k\geq 1} \frac{m / k^p}{(2m-1)^2 - 4k^2}.
\end{align}
Now, for any given \(m,p\) one can compute \(E_{m,p}\) to high precision using
Euler-Maclaurin summation. Then for fixed \(\vec{a}\), the sum over \(m\) can
be estimated to high precision using series acceleration methods such as the
Levin transform (see, for example, \cite{sidi2003}). In practice, these means
that one needs to pre-compute a set of \(E_{m,p}\) for \(p=0,\dots,P\) and
\(m=1,2,\dots,N\) precise to many digits. We typically used something like
\(P=32, N=4096\) with computations to 1024 digits --- not a small computation,
but by no means an onerous one on modern hardware.

So, for a given fixed vector \(\vec{a}\) the above expressions allow us to
estimate \(\mathcal{C}(\vec{a})\) to high precision. We note that by
differentiating the above expressions, one can similarly estimate the partial
derivatives \(\pdiff{}{a_j}\mathcal{C}\), giving the gradient, and (in
principle) the  Hessian. This means that one can employ numerical methods to
optimise over \(\vec{a}\). We did this using the BFGS quasi-Newton method (see, say, \cite{fletcher2013}). We
made these computations using the \texttt{mpmath} package for \texttt{python}~\cite{mpmath}.
That library has many arbitrary precision numerical methods, rendering much of
the computation quite easy to code.

This lead us to very precise estimates of a near-optimal vector of coefficients
\(\vec{a}\), and so a near-optimal \(f\)-function. We found that the estimated
\(a_j\) were quite stable under increasing \(P\). We also examined the ansatz
\begin{equation}
	\hat{f}(k) = (-1)^k \left(
	\frac{1}{\sqrt{k}} \sum_{i=0}^{P-1} a_i k^{-i}
	+ \frac{1}{k} \sum_{i=0}^{P-1} b_i k^{-i}
	\right)
\end{equation}
giving both integer and half-integer powers of \(k\). However, in that case
we found that the \(b_i\) in the above expression were extremely close to zero,
especially as \(P\) was increased.

We struggled for some time with the problem of taking a fixed (hopefully near-optimal) \(\vec{a}\)
and turning that into a rigorously bounded value for \(\mathcal{C}(\vec{a})\). Unfortunately
the double-sums involving the \(E_{m,p}\) made this rather challenging (to say the least). Eventually
we gave up on this approach. However, examining the near-optimal \(f\)-function led us to a simpler ansatz.

\begin{figure}
	\begin{center}
		\includegraphics[width=0.45\textwidth]{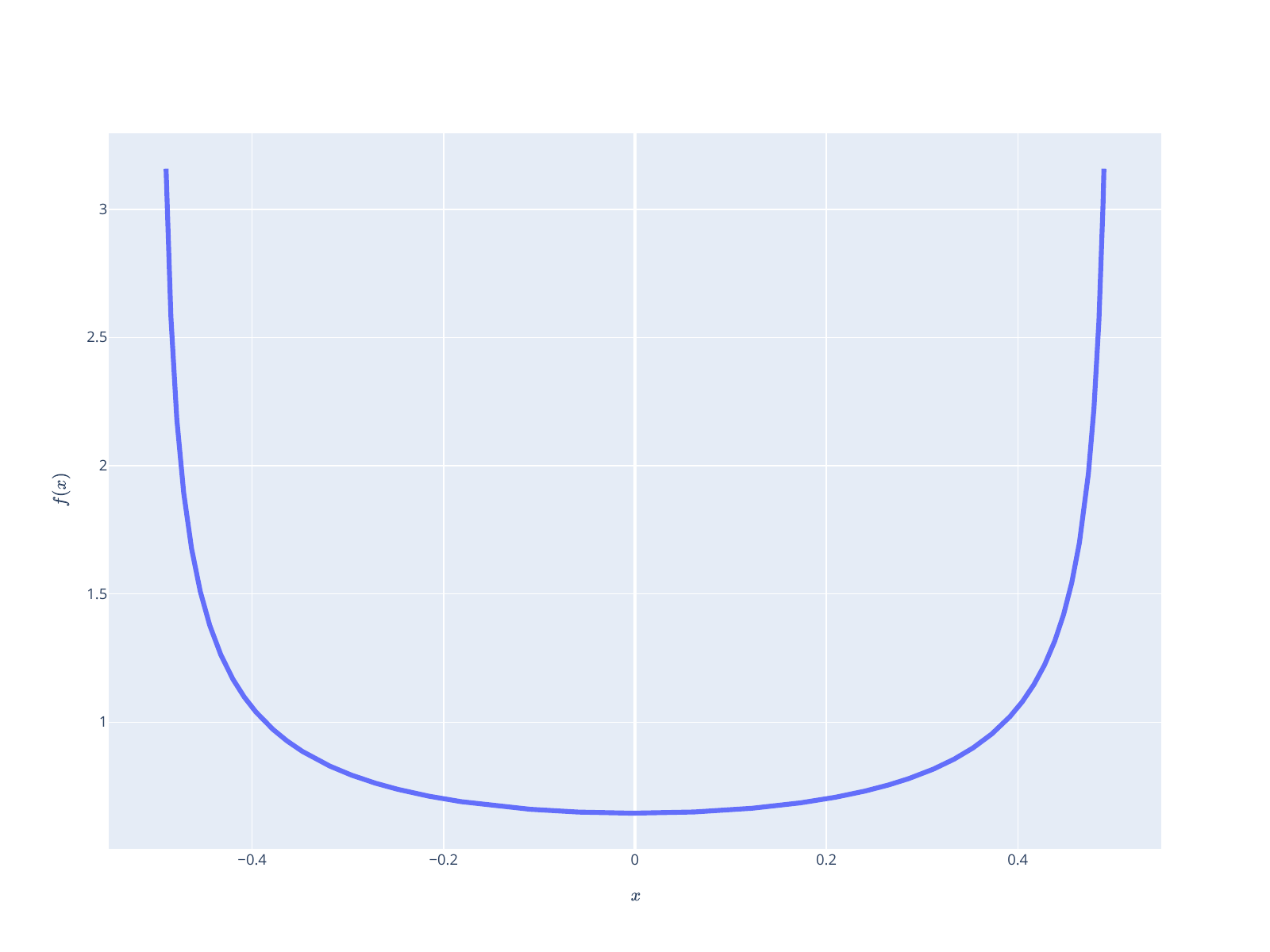}
		\includegraphics[width=0.45\textwidth]{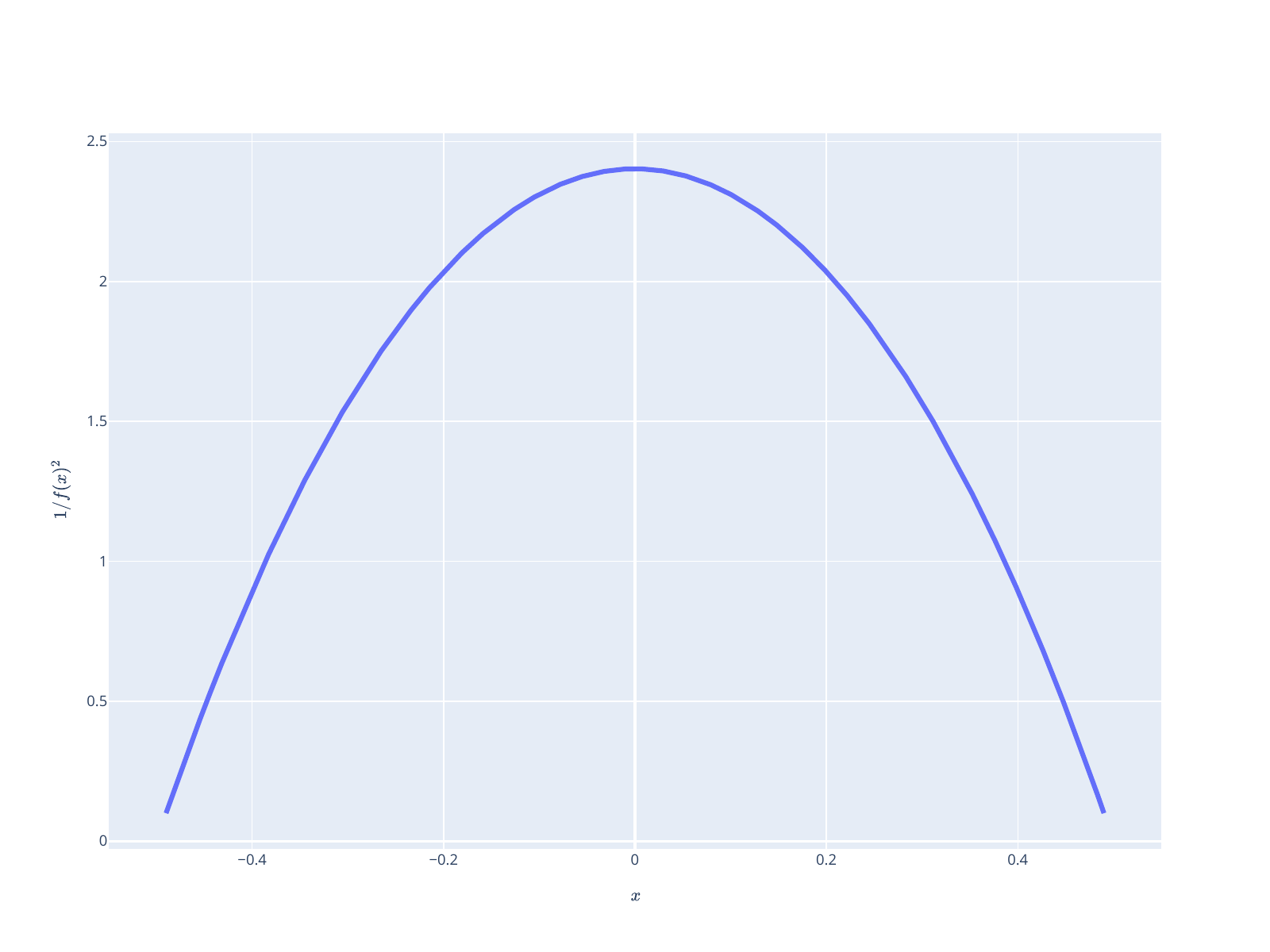}
	\end{center}
	\caption{The near-optimal \(f(x)\) as computed via the ansatz in equation~\eqref{eqn_ansatz_1}. We also
		plot \(1/f(x)^2\) to show the nature of the singularity as \(x \to \pm 1/2\). This second plot
		is strongly suggestive that \(1/f(x)^2 \approx \text{const}\cdot(1-4x^2)\).}
	\label{fig_near_opt}
\end{figure}

From the \(\vec{a}\) coefficients, we can compute the corresponding \(f(x)\)
\begin{align}
	f(x) & = \sum_{k \in \ZZ} \hat{f}(k) \cdot e^{2\pi i k x}
	= 1 + 2 \sum_{j=0}^{P-1} a_j \sum_{k\geq 1} \frac{(-1)^k}{k^{j+1/2}} \cos(2 \pi k x) \\
	     & = 1 + 2 \sum_{j=0}^{P-1} a_j \claus(j+1/2, 2 \pi x + \pi )
	\qquad \text{ with }
	\claus(p,\theta)  = \sum_{k \geq 1} \frac{\cos k\theta}{k^p},
\end{align}
where \(\claus(p,\theta)\) is a Clausen cosine function.
The \texttt{mpmath} package then allows us to estimate the near optimal \(f(x)\) at points in its domain,
which we plot in Figure~\ref{fig_near_opt}.

\begin{figure}
	\begin{center}
		\includegraphics[width=0.45\textwidth]{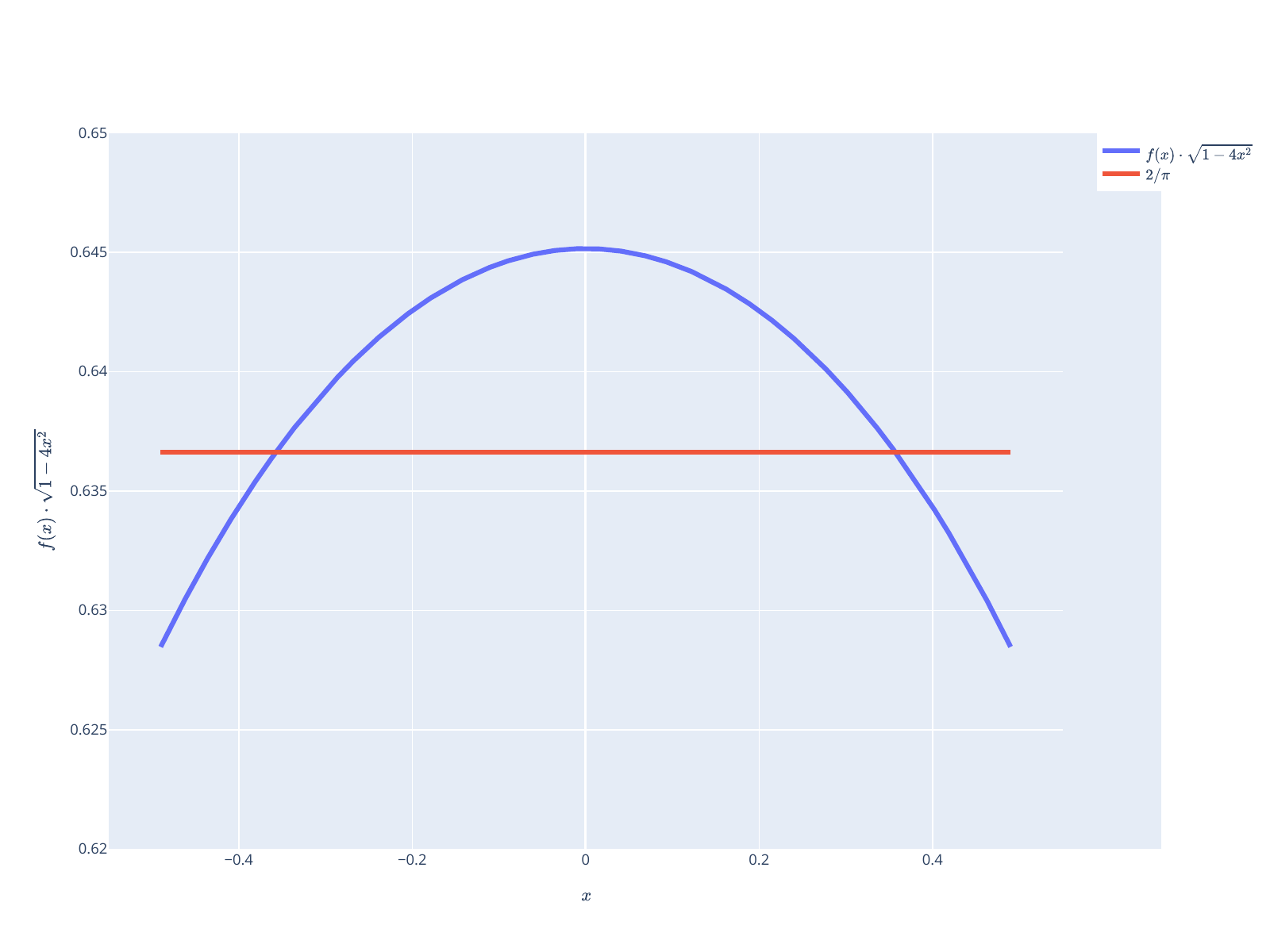}
		\includegraphics[width=0.45\textwidth]{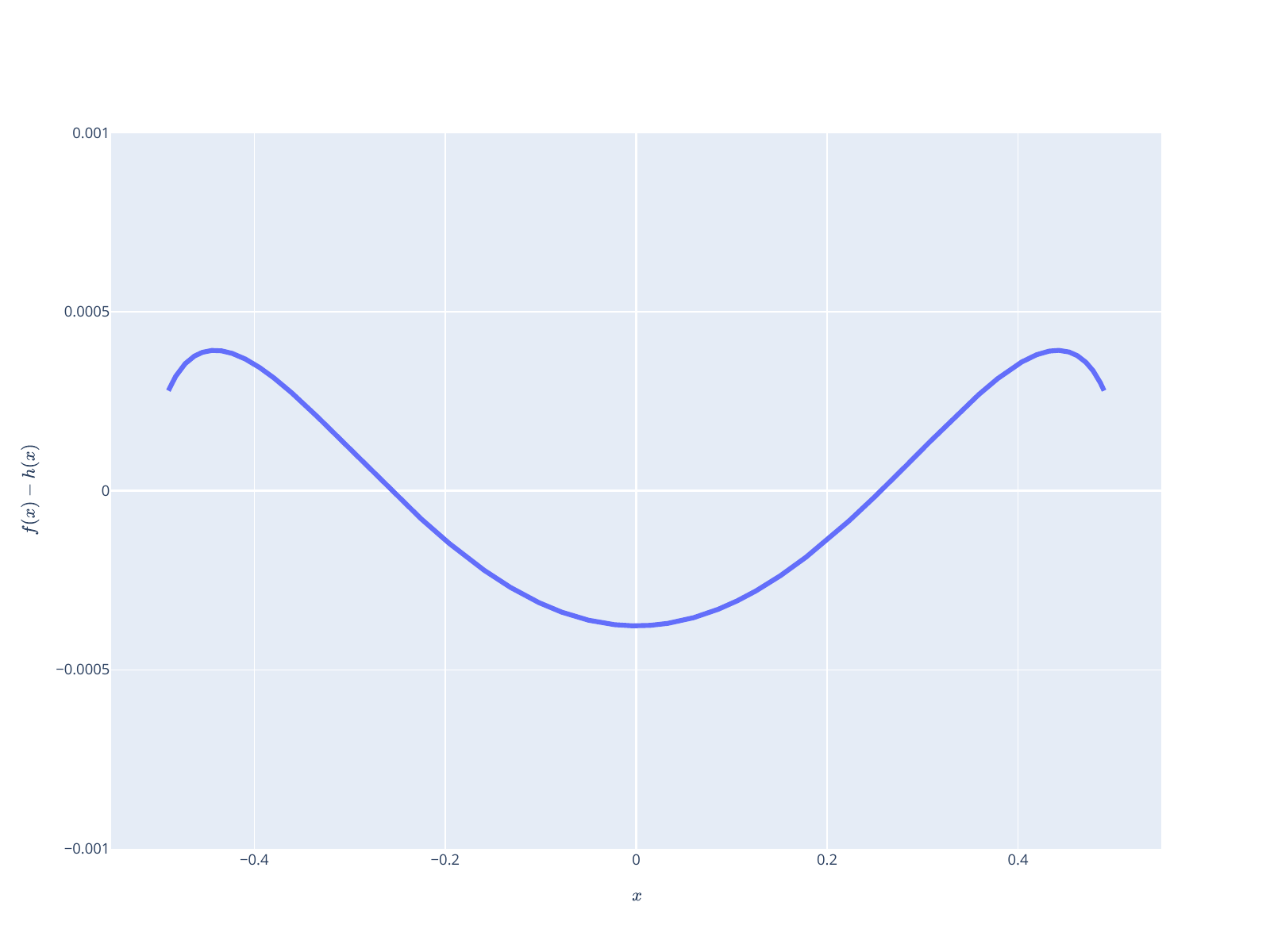}
	\end{center}
	\caption{A plot of \(f(x) \cdot \sqrt{1-4x^2}\) with a straight line at \(y=2/\pi\); note the vertical range is \([0.62,0.65]\).
		In the second figure we plot \(f(x)-h(x)\) with \(h(x)\) given by equation~\eqref{eqn_approx_f_2} and find that it
		is very close to \(0\); note the vertical range is \([-0.001,0.001]\).}
	\label{fig_near_opt_ansatz}
\end{figure}

This function quite clearly diverges to \(+\infty\) as \(x \to \pm 1/2\). We see strong evidence that this
is an inverse square-root singularity by also plotting \(1/f(x)^2\); the resulting curve looks very much as though it is
proportional to the parabola \((1-4x^2)\). Indeed, this suggests that the near-optimal \(f\) should be approximately
\begin{equation}
	h(x) = \frac{1}{\sqrt{1-4x^2}} \cdot \frac{2}{\pi}
	\label{eqn_approx_f_1}
\end{equation}
where the multiplicative constant ensures that the function has mass \(1\). We plot \(f(x)\cdot \sqrt{1-4x^2}\)
in Figure~\ref{fig_near_opt_ansatz},
and see that the resulting curve is not too far from the straight line \(y=2/\pi\).
We note that this was observed by White~\cite{white_data, white2022}, who also credits O'Bryant for suggesting
the family of functions \( (1-4x^2)^c \). This observation helped prompt our explorations.

Some simple curve-fitting
shows that \(f(x) \approx h(x)\) given by
\begin{equation}
	h(x) = \frac{2}{\pi} \frac{1}{\sqrt{1-4x^2}} \cdot 0.986 + \frac{4}{\pi} \sqrt{1-4x^2} \cdot 0.014
	\label{eqn_approx_f_2}
\end{equation}
Again, the multiplicative constants are chosen so \(\int h = 1\).
In Figure~\ref{fig_near_opt_ansatz} we also plot \(f(x) - h(x)\) and find that it is very close to \(0\).

Note that if we set \(h(x) = 2/\pi \cdot (1-4x^2)^{-1/2}\) (as per \eqref{eqn_approx_f_1}), then following equation~\eqref{eqn_big_f_hat} we
obtain
\begin{equation}
	\hat{H}(k) = \frac{1}{2} \bess{0}{\frac{\pi k}{2}}
\end{equation}
where \(\bess{p}{z}\) is the Bessel function of first kind of order \(p\). Then by equation~\eqref{eqn_c_by_big_f}, we have
\begin{equation}
	\mathcal{C}(h) = \frac{1}{2} + \sum_{k \geq 1} \bess{0}{\frac{\pi k}{2} }^4 \approx 0.574694862\dots
\end{equation}
This number constitutes an upper bound on \(\nu_2^2\) and is only slightly
above White's bound.

We note that the number quoted above is an approximation of the sum and not entirely rigorous. Below we will discuss how
we can compute rigorous error-bars in such summations.  We also note that although the summands are all positive, they have
a strong period-4 behaviour (more on this below), and we found it advantageous to rewrite the sum as
\begin{equation}
	\sum_{k \geq 1} \bess{0}{\frac{\pi k}{2} }^4  =	\sum_{n \geq 0} \left( \sum_{j=1}^4  \bess{0}{\frac{\pi(4n+j) }{2}}^4 \right)
\end{equation}
when doing numerical explorations --- grouping this way makes the outer partial sums monotonic.

Now if we set \(h(x)\) by equation~\eqref{eqn_approx_f_2}, then we have
\begin{equation}
	\hat{H}(k) = \frac{1}{2} \left( \bess{0}{\frac{\pi k}{2}} \cdot 0.986 + \frac{4}{\pi k} \cdot \bess{1}{\frac{\pi k}{2}} \cdot 0.014 \right).
	\label{eqn_approx_fhat_2}
\end{equation}
Again, using~\eqref{eqn_c_by_big_f}, we have
\begin{equation}
	\mathcal{C}(h) = \frac{1}{2} + \sum_{k \geq 1}
	\left( \bess{0}{\frac{\pi k}{2}} \cdot 0.986 + \frac{4}{\pi k} \cdot \bess{1}{\frac{\pi k}{2}} \cdot 0.014
	\right)^4
	\approx 0.5746396478\dots
\end{equation}
This is already inside the bounds computed by White and, if the above summation is made rigorous, would constitute
a tighter upper bound.

\section{A second ansatz and rigorous upper bounds}\label{sec second ansatz}
We can improve the above further by setting
\begin{equation}
	\label{eqn_2nd_ansatz_start}
	h(x) = \frac{2}{\pi} \frac{1}{\sqrt{1-4x^2}} \cdot a + \frac{4}{\pi} \sqrt{1-4x^2} \cdot (1-a),
\end{equation}
which then gives
\begin{align}
	\mathcal{C}(h;a) & = \frac{1}{2} + \sum_{k \geq 1}
	\left( \bess{0}{\frac{\pi k}{2}} \cdot a + \frac{4}{\pi k} \cdot \bess{1}{\frac{\pi k}{2}} \cdot (1-a)
	\right)^4                                                                                    \\
	                 & \approx 0.07585490614 a^4 - 0.1340881050 a^3 + 0.2466255031 a^2 \nonumber \\
	                 & \quad - 0.3863025578 a  + 0.7726051152. \nonumber
\end{align}
This quartic in \(a\) is minimised when \(a \approx 0.9863131909 \) which gives
\begin{equation}
	\mathcal{C}(h)  \approx 0.5746396188
	\label{eqn_approx_upper}
\end{equation}
All the numbers in the above expressions can be computed to higher precision.
We can improve the result much further by extending~\eqref{eqn_2nd_ansatz_start} to the
ansatz
\begin{equation}
	f(x) =  \sum_{j=0}^{P-1} a_j \cdot \binom{j}{1/2} \cdot (1-4x^2)^{j-1/2} \qquad 1=\sum a_j
	\label{eqn_ansatz_2}
\end{equation}
where the multiplicative constants and the additional constraint ensure that \(\int f =1\). This, in turn
gives an ansatz for the Fourier coefficients
\begin{equation}
	\hat{F}(k) = \frac{1}{2} \sum_{j=0}^{P-1} a_j \cdot \bess{j}{\frac{\pi k}{2}} \cdot j! \left( \frac{4}{\pi k} \right)^j.
\end{equation}
Accordingly we define
\begin{equation}
	\mybess{p}{k} =  \bess{p}{\frac{\pi k}{2}} \cdot p! \left( \frac{4}{\pi k} \right)^p,
\end{equation}
and so we can write
\begin{equation}
	\mathcal{C}(\vec{a})  = \frac{1}{2} + \sum_{k \geq 1}
	\left(
	\sum_{j=0}^{P-1} a_j \cdot \mybess{j}{\frac{\pi k}{2}}
	\right)^4.
	\label{eqn_ansatz_2k}
\end{equation}
Our problem can now be restated as one of minimising the above multinomial with \(\binom{P+4}{4}\) terms,
subject to the constraint \(\sum a_j =1\).

For small values of \(P\) it is feasible to pre-compute the \(\binom{P+4}{4}\) sums of products of Bessel functions,
and use Newton's method with Lagrange multipliers to find a near-optimal \(\vec{a}\) and an upper
bound on \(\nu_2^2\) --- of course one still needs to bound any floating point errors in the final result.
We did try this approach, but found that the quartic growth in multinomial terms became too prohibitive at
modest values of \(P\).

Instead we worked directly with equation~\eqref{eqn_ansatz_2k}. This shifts the problem to
one of computing this sum (and its partial derivatives) extremely efficiently at any given fixed value of \(\vec{a}\).
To do this we use Kummer's series transform (see \cite{abramowitz1964}~p.16) and the asymptotics of Bessel functions. So start by writing our sum as
\(\mathcal{C} = \sum s_n\) and then we find a similar, but summable series \(B = \sum b_n\), so that
the \(b_n\) has the same dominant asymptotics as \(s_n\). Then we can compute \(A = B + \sum_n (s_n-b_n)\).
Since \(s_n \sim b_n\), the summands \((s_n-b_n)\) decay to zero faster and so the series converges more quickly.
We can bounds \(\abs{s_n-b_n}\) by terms of the form \(\frac{const}{n^K}\). To ensure that these terms are small,
we will compute
\begin{equation}
	\sum_{n=1}^\infty s_n = \sum_{n=1}^{N} s_n + \sum_{n=N+1}^\infty b_n
	+ \sum_{n=N+1}^\infty (s_n-b_n)
	\label{eqn_kummer_transform}
\end{equation}
and then bound
\begin{equation}
	\abs{\sum_{n=N+1}^\infty (s_n-b_n)} \leq \const \cdot \sum_{n=N+1}^\infty n^{-K} = \const \cdot \zeta(K, N)
	\leq \const \cdot N^{1-K}.
	\label{eqn_kummer_error}
\end{equation}

In practice the above idea is complicated by the periodic behaviour of the asymptotics of Bessel functions.
We use a standard large-argument expansion
\begin{align}
	\mybess{p}{\frac{\pi k}{2}}
	            & \sim \frac{2^{2p+1}}{\pi^{p+1}} \cdot \frac{p!}{k^{p+1/2}} \cdot \left( \cos \omega \sum_{\ell\geq0} (-1)^\ell \frac{v_{2\ell}(p)}{k^{2\ell}}
	- \sin\omega \sum_{\ell\geq 0} (-1)^\ell \frac{v_{2\ell+1}(p)}{k^{2\ell+1}} \right)                                                                         \\
	v_{\ell}(p) & = \frac{(-1)^p}{4^\ell \pi^\ell \ell!} \prod_{j=1}^\ell (4p^2-(2j-1)^2)
	\quad \text{and} \quad \omega   = \frac{\pi (k-p) }{2} - \frac{\pi}{4}
\end{align}
where the error in this expansion is bounded by the magnitude of the first neglected term in both sums,
provided the index of that term is at least \(p/2\) and \(k \geq 1\) (see, say, \cite{olver1997} or \cite{olver2019} (10.17.i)). For integer \(k\), this expansion has a period-4 sign pattern, and so in practice we compute 4 separate expansions:
\begin{align}
	\sum_{j=0}^{P-1} a_j \cdot \mybess{j}{\frac{\pi k}{2}}
	 & \sim \frac{\sqrt{2}}{\pi \sqrt{k}} \cdot \begin{cases}
		                                            +a_0 - \frac{a_0+16a_1}{4\pi k} - O(k^{-2}) & k \equiv 0 \\
		                                            +a_0 + \frac{a_0+16a_1}{4\pi k} - O(k^{-2}) & k \equiv 1 \\
		                                            -a_0 + \frac{a_0+16a_1}{4\pi k} + O(k^{-2}) & k \equiv 2 \\
		                                            -a_0 - \frac{a_0+16a_1}{4\pi k} + O(k^{-2}) & k \equiv 3
	                                            \end{cases}.
	\label{eqn_bess_res}
\end{align}
From those, we compute their 4th powers:
\begin{align}
	\left(\sum_{j=0}^{P-1} a_j \cdot \mybess{j}{\frac{\pi k}{2}} \right)^4
	 & \sim \frac{4}{\pi^4 k^2} \cdot \begin{cases}
		                                  +a_0^4 - \frac{a_0^3(a_0+16a_1)}{\pi k} - O(k^{-2}) & k \equiv 0 \\
		                                  +a_0^4 + \frac{a_0^3(a_0+16a_1)}{\pi k} - O(k^{-2}) & k \equiv 1 \\
		                                  +a_0^4 - \frac{a_0^3(a_0+16a_1)}{\pi k} - O(k^{-2}) & k \equiv 2 \\
		                                  +a_0^4 + \frac{a_0^3(a_0+16a_1)}{\pi k} - O(k^{-2}) & k \equiv 3
	                                  \end{cases}
	\label{eqn_bess4_res}
\end{align}
We note that the \(O(\cdot)\) term is rigorously bounded as a term of the form \(\frac{cosnt}{k^2}\).
More generally we terminate these expansions at \(k^{-K}\).

Since we have a bound on the \(O(\cdot)\)  we can sum the expansion in four parts and keep track of the error term.
The sum of the four expansions splits into sums of the form
\begin{equation}
	\sum_{n=N}^\infty \frac{1}{(4n+j)^s}
	= 4^{-s} \sum_{n \geq 0} \frac{1}{(n+N+j/4)^s} = 4^{-s} \zeta(s,N+j/4)
	\label{eqn_hurwitz_four}
\end{equation}
and so expressed in terms of the Hurwitz zeta-function. So, for a given fixed \(\vec{a}\) we
\begin{itemize}
	\item Compute the four expansions in equation~\eqref{eqn_bess_res}.
	\item From those compute the four expansions in equation~\eqref{eqn_bess4_res}.
	\item Sum the individual terms over \(k\gt N\) in the expansion using equation~\eqref{eqn_hurwitz_four}.
	\item Similarly bound the error arising from truncating the asymptotic expansions as per equation~\eqref{eqn_kummer_error}.
	\item Compute a finite sum of equation~\eqref{eqn_ansatz_2k} for \(1 \leq k \leq N\).
	\item Add all of these contributions.
\end{itemize}
and so arrive at \(\mathcal{C}(\vec{a})\) with a bound on the error.

To control floating point errors we coded all of the above in \texttt{c++} and made heavy use of the \texttt{flint} library~\cite{flint}.
This library has arbitrary precision floating point routines using ball-arithmetic~\cite{johansson2017}. This facilitates mathematically rigorous
floating point computations, and allows us compute \(\mathcal{C}(\vec{a})\) with rigorous bounds on any floating point error.

To now minimise over \(\vec{a}\) we use Lagrange multipliers
\begin{equation}
	\mathcal{L}(\vec{a}, \lambda) = \mathcal{C}(\vec{a})	+ \lambda \left(1 - \sum_p a_p \right)
	\label{eqn_lagrangian}
\end{equation}
and Newton-Raphson to solve \(\nabla \mathcal{L} = 0\). This requires us to compute
the first and second partial derivatives of \(\mathcal{L}(\vec{a})\):
\begin{align}
	\pdiff{}{a_\ell} \mathcal{L}(\vec{a},\lambda)
	                                               & =
	4 \sum_{k \geq 1} \mybess{\ell}{\frac{\pi k}{2}}
	\left(
	\sum_{j=0}^{P-1} a_j \cdot \mybess{j}{\frac{\pi k}{2}}
	\right)^3 - \lambda, \nonumber                                   \\
	\pdiff{}{a_i}\pdiff{}{a_\ell} \mathcal{L}(\vec{a},\lambda)
	                                               & =
	12 \sum_{k \geq 1} \mybess{i}{\frac{\pi k}{2}}\mybess{\ell}{\frac{\pi k}{2}}
	\left(
	\sum_{j=0}^{P-1} a_j \cdot \mybess{j}{\frac{\pi k}{2}}
	\right)^2, \nonumber                                             \\
	\pdiff{}{\lambda} \mathcal{L}(\vec{a},\lambda) & = - \sum_p a_p,
	\qquad
	\pdiff{}{\lambda} \pdiff{}{a_\ell} \mathcal{L}(\vec{a},\lambda)     = -1,
	\qquad
	\frac{\partial^2}{\partial \lambda^2} \mathcal{L}(\vec{a},\lambda)  = 0.
	\label{eqn_lagrange_multipliers}
\end{align}
We can compute similar expansions to those in equation~\eqref{eqn_bess4_res}, and then
compute the sum via much the same procedure used to compute \(\mathcal{C}(\vec{a})\). From these quantities we
can then perform Newton-Raphson iterations to converge towards the optimal \(\vec{a}\)-value.

Since many \(\zeta\)-function values are reused in each iteration, as well as
the expansions of \(\mybess{j}{\frac{\pi k}{2}}\), we pre-compute them. This
speeds up the repeated computations considerably. We parallelise the computation
of the partial derivatives using a thread-pool provided by the \texttt{flint} library.
All of this means that the bottle-neck of the computation is solving for the step direction \(\vec{s}\) via
\begin{equation}
	\mathcal{H} = \nabla^2 f  \qquad \vec{s} = - \mathcal{H}^{-1} \nabla f.
\end{equation}
Using this approach we computed near-optimal \(\vec{a} = (a_0, \dots, a_{P-1}), \lambda\), and the corresponding \(\mathcal{C}(\vec{a})\),
for a range of different \(P\). Each of these constitutes an upper bound for \(\nu_2^2\).

\section{Lower bound} \label{sec lower}
Following White (see \cite{white2024} Lemma~3.2), we show how a precise upper bound can be leveraged to construct a similarly precise lower bound.
Define a new 1-periodic function \(g:(-1/2,1/2) \to \RR\) so that \(\int g =2\) and then extend it to \(G:[-1,1] \to \RR\)  by
\begin{equation}
	G(x) = \begin{cases}
		1      & x \in [-1/2,1/2]  \\
		1-g(x) & \text{ otherwise}
	\end{cases}.
	\label{eqn_lower_bound_G}
\end{equation}
This means that
\begin{align}
	\int_{-1}^1 F(x) G(x) \dee{x} & = \int_{-1/2}^{1/2} f(x) \dee{x} = 1.
	\label{eqn_int_of_FG}
\end{align}
Observe that \(\hat{G}(0)=0\) and then Plancherel's theorem then implies that
\begin{align}
	1 = \int_{-1}^1 F(x) \overline{G(x)} \dee{x}
	 & = 2 \sum_{k\neq 0} \hat{F}(k) \overline{\hat{G}(k)}.
\end{align}
Applying H\"older's inequality we get
\begin{equation}
	\frac{1}{16} \leq \left( \sum_{k \neq 0} \abs{\hat{F(k)}}^4 \right) \left( \sum_{k \neq 0} \abs{\hat{G}(k)}^{4/3} \right)^3.
\end{equation}
Rearranging this inequality and substituting the optimal \(F\) we have
\begin{equation}
	\nu_2^2 \geq \frac{1}{2} + \frac{1}{2} \left( \sum_{k \neq 0} \abs{\hat{G}(k)}^{4/3} \right)^{-3}.
	\label{eqn_nu_lower_bound}
\end{equation}
Consequently to form a good lower bound we need to optimise \(\hat{G}\), so that H\"older's inequality is as close to an equality
as possible. It becomes an equality when \(\hat{G}(k) = \alpha \hat{F}(k)^3\) and so we seek a function \(G(x)\) whose
Fourier coefficients are proportional to the cube of those of \(F(x)\).

We note that \(\hat{F}(k)^3\) is the Fourier transform of \((F \ast F \ast F)(x)\), which is a function that
arises naturally when proving the uniqueness of the optimising function (see Lemma~28 in~\cite{green2001} and Lemma~6 in~\cite{white2024}). For the optimal \(F(x)\), the function \((F \ast F \ast F)(x)\) is
constant on \((-1/2,1/2)\) with value \(\nu_2^2 / 4\). That analysis also shows that the correct value of
\(\alpha\) in the above equation is \(\alpha = 8/(2\nu_2 -1) \approx 53.8\).

\begin{figure}
	\begin{center}
		\includegraphics[width=0.45\textwidth]{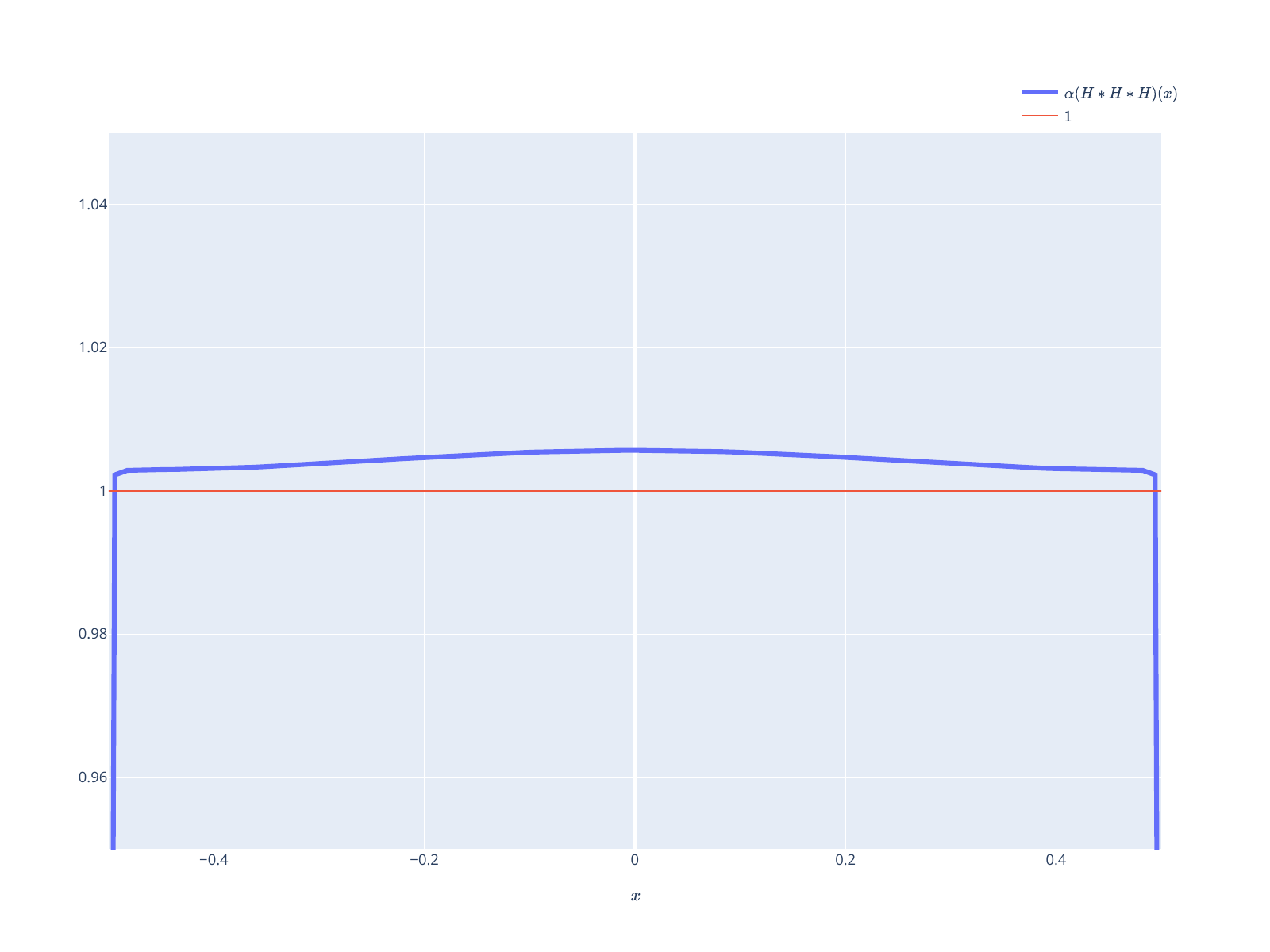}
		\includegraphics[width=0.45\textwidth]{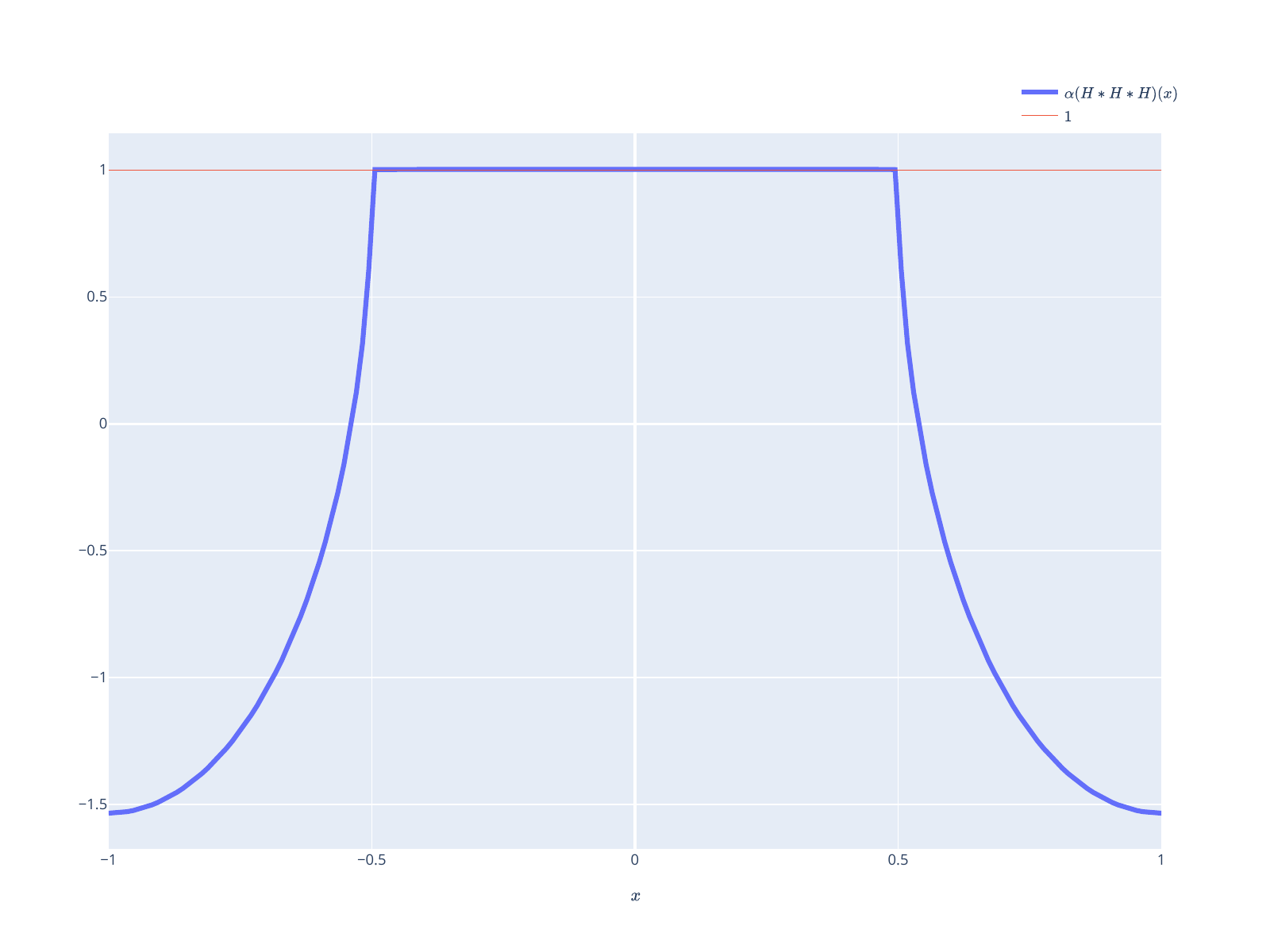}
	\end{center}
	\caption{A plot of \((H \ast H \ast H)(x)\) based on \(H(x)\) from equation~\eqref{eqn_approx_f_2} for \( x \in (-1/2,1/2) \).
		Notice that the function is fairly constant on this range and takes a value close to \(1\).
		The second plot shows the same function but for \(x \in (-1,1)\).}
	\label{fig_near_opt_H_cube}
\end{figure}

When we take a near-optimal \(F\) we observe the result is nearly constant. For example, using equations~\eqref{eqn_approx_f_2} and~\eqref{eqn_approx_fhat_2},
we compute the Fourier series of \((H \ast H \ast H)(x)\) as the cube of \(\hat{H}(k)\) and multiply by \(\alpha=53.8\):
\begin{equation}
	\alpha \hat{H}(k)^3 = \frac{53.8}{8} \left( \mybess{0}{\frac{\pi k}{2}} \cdot 0.986 + \frac{4}{\pi k} \cdot \mybess{1}{\frac{\pi k}{2}} \cdot 0.014 \right)^3.
\end{equation}
Unfortunately these coefficients do not easily transform back to a simple function of \(x\), however we can use it to construct an approximate plot of \(\alpha (H \ast H \ast H)(x)\).
We compute the first few hundred Fourier coefficients and then use those to construct a truncated Fourier-cosine series for \( (H \ast H \ast H)(x)\)
and then plot that in Figure~\ref{fig_near_opt_H_cube}.
That quite clearly shows the function is nearly constant on \((-1/2,1/2)\) taking a value very nearly \(\nu_2^2\).

This indicates a problem, but does also help us towards a solution. If we have exactly the optimal \(F\),
then we can set \(\hat{G} \propto \hat{F}^3\) and obtain a bound via H\"older's inequality as above. However, we only have a near-optimal \(F\),
and if we use that to construct a \(G\), then it will not satisfy the requirements of equation~\eqref{eqn_lower_bound_G}.
We cannot guarantee that the \(G(x)\) will be constant on \((-1/2,1/2)\) and so we cannot assert that \(\int F G =1\) as we require in equation~\eqref{eqn_int_of_FG}.
Additionally, the correct constant \(\alpha\) requires a knowledge of \(\nu_2^2\) which is precisely the constant we are trying to compute!

\begin{figure}
	\begin{center}
		\includegraphics[width=0.45\textwidth]{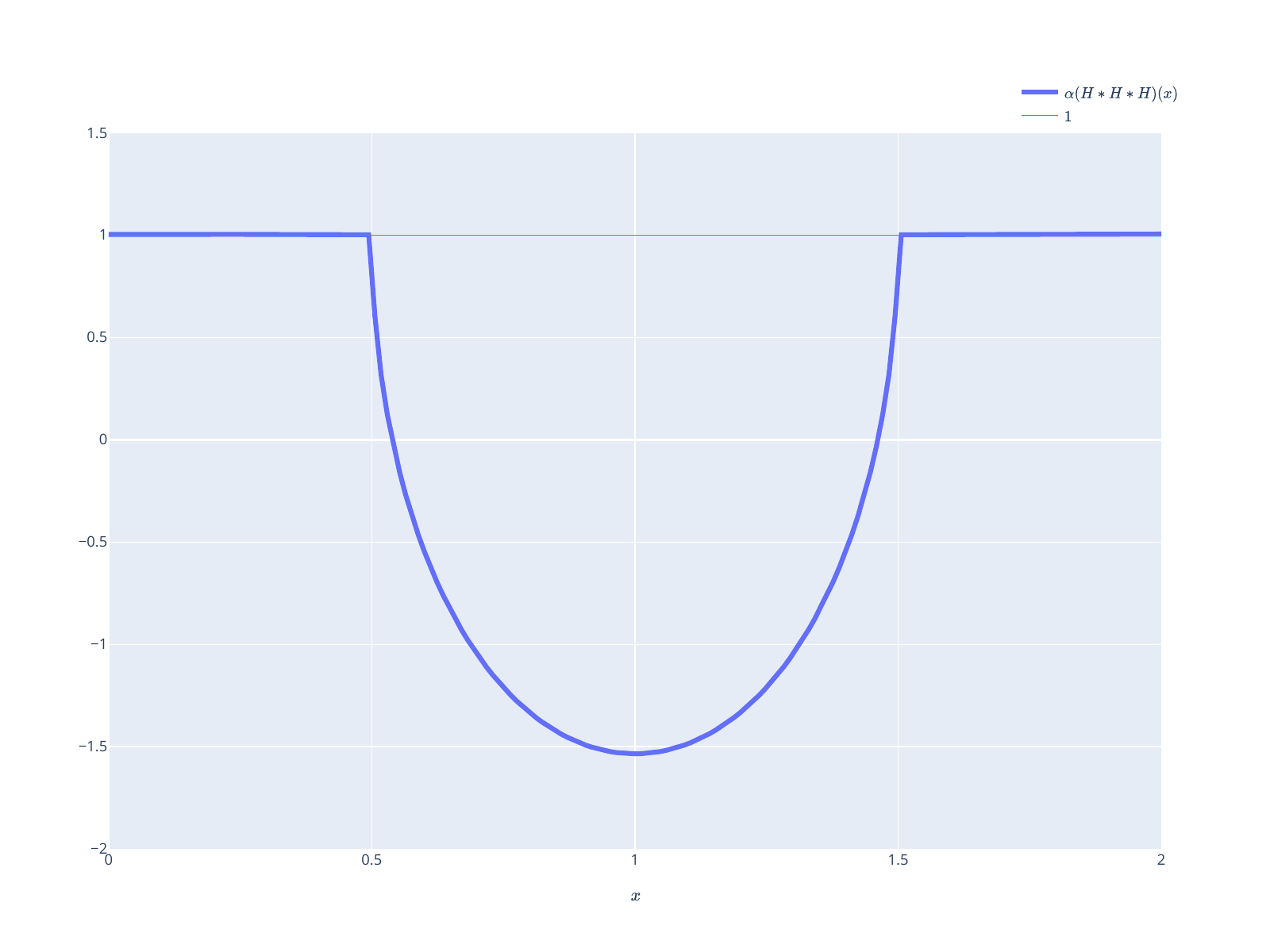}
		\includegraphics[width=0.45\textwidth]{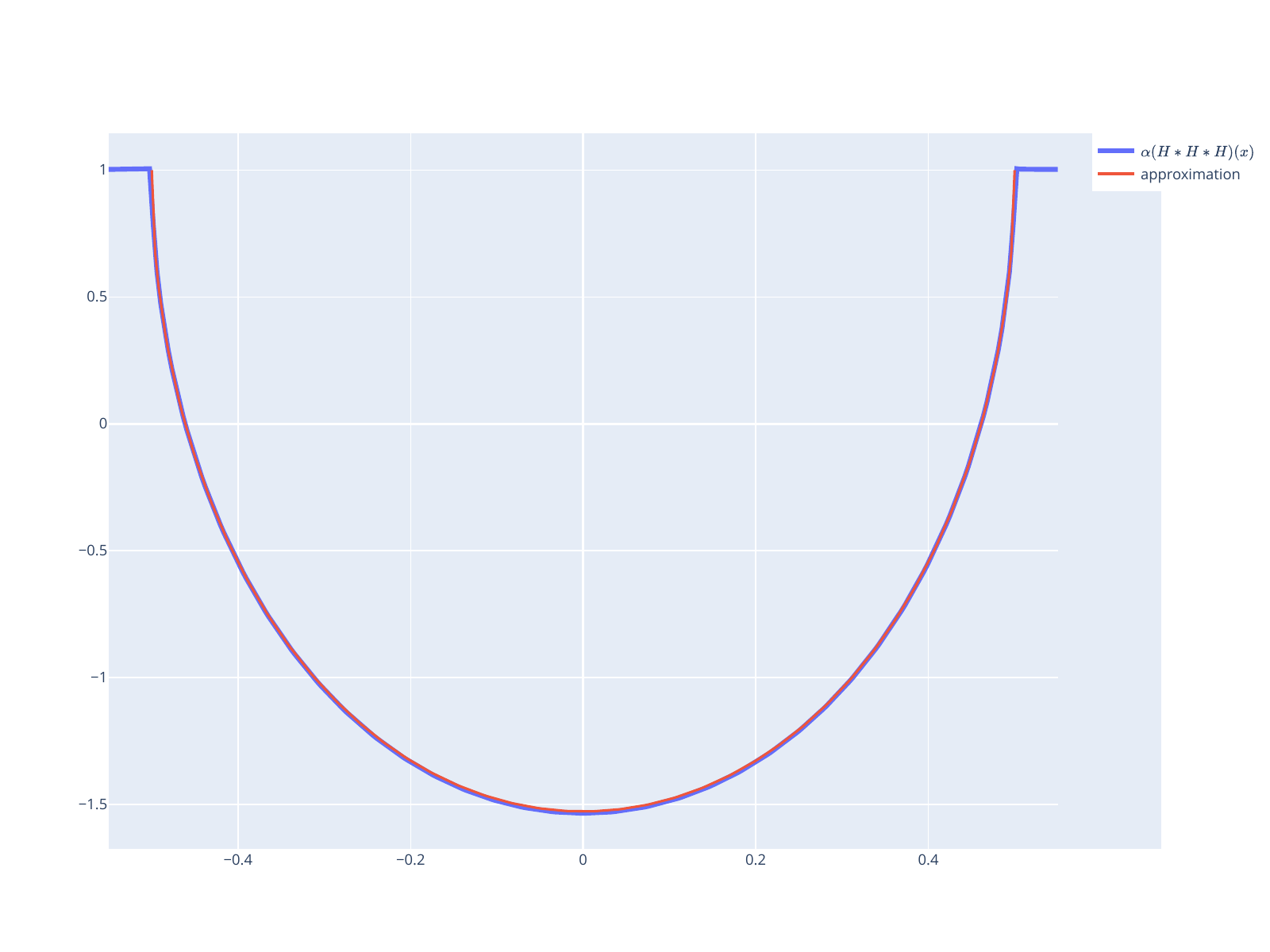}
	\end{center}
	\caption{A plot of \((H \ast H \ast H)(x)\)from equation~\eqref{eqn_approx_f_2}, but now
		shifted by 1 to highlight the non-constant region. The second plot shows that this function is well approximated by equation~\eqref{eqn_h_cube_approx}. }
	\label{fig_near_opt_H_cube2}
\end{figure}

In Figure~\ref{fig_near_opt_H_cube2} we shift the plot by 1 to highlight
the non-constant region. Some simple least-squares curve fitting shows that the non-constant region is well approximated by a
function of the form \(1 + b_1 \sqrt{1-4x^2} + b_2 (1-4x^2)^{3/2}\). In particular
\begin{equation}
	\alpha (H \ast H \ast H)(x+1) \approx 1 - 2.044 \cdot \frac{4}{\pi} (1-4x^2)^{1/2} + 0.044 \cdot \frac{16}{3\pi} (1-4x^2)^{3/2},
	\label{eqn_h_cube_approx}
\end{equation}
where we have used a similar expansion to that used in equation~\eqref{eqn_approx_f_2}. Notice that the right-hand side
of the above equation integrates to -1, as is required for equation~\eqref{eqn_lower_bound_G}. Additionally, the Fourier transform of the right-hand side is easily expressed in terms of
\(\mybess{1}{\frac{\pi k}{2}}, \mybess{2}{\frac{\pi k}{2}}\). Then notice that for large \(k\) we have
\begin{equation}
	(-1)^k \left(	2.044 \cdot\frac{1}{2}  \mybess{1}{\frac{\pi k}{2}} - 0.044 \cdot\frac{1}{2}  \mybess{2}{\frac{\pi k}{2}} \right)
	\approx \alpha \hat{H}(k)^3,
\end{equation}
where the factor of \((-1)^k\) accounts for the shift \(x \mapsto x+1\). Summing the above (non-rigorously) gives
\begin{equation}
	\sum_{k \neq 0} \abs{1.022 \cdot \mybess{1}{\frac{\pi k}{2}} - 0.022 \cdot \mybess{2}{\frac{\pi k}{2}} }^{4/3}
	\approx 1.885096418\dots
	\qquad
	\nu_2^2 \geq 0.5746395994\dots
\end{equation}
Combining this with equation~\eqref{eqn_approx_upper} gives
\begin{equation}
	0.5746395995\dots \leq \nu_2^2 \leq 0.5746396188\dots
\end{equation}
which are tighter than the bounds computed by White, but, at this stage we are yet to control for summation and floating-point errors.

We use this observation to formulate an ansatz for \(G(x)\).  We construct
\begin{equation}
	g(x) = \sum_{p=1}^{P} b_p (1-4x^2)^{p-1/2} \cdot \binom{p}{1/2}.
\end{equation}
We then satisfy \(\int g=2\) by requiring that \(\sum b_p = 2\). The Fourier coefficients of the corresponding \(G\) are then
\begin{equation}
	\hat{G}(k) = (-1)^k \sum_{p=1}^{P} b_p\cdot \frac{1}{2}\mybess{p}{\frac{\pi k}{2}},
\end{equation}
where the factor of \((-1)^k\) accounts for the shift required by equation~\eqref{eqn_lower_bound_G}.
In order to determine the coefficients \(b_p\), we require that the asymptotic expansion of \(\hat{G}(k)\) agrees
with the expansion of \(\hat{F}(k)^3\).

More precisely we take our near-optimal \(\vec{a}\), and then compute an expansion
in equation~\eqref{eqn_bess_res} and its cube. To leading order this is \(O(n^{-3/2})\), and
so can be matched to the leading order term of \(\mybess{1}{\frac{\pi k}{2}}\) to compute the coefficient
\(b_1\). We then determine \(b_2\) from the coefficient of \(n^{-5/2}\) which is a combination of
the leading term of \(\mybess{2}{\frac{\pi k}{2}}\) and the next-to-leading term of \(\mybess{1}{\frac{\pi k}{2}}\).
We continue to the coefficient of \(n^{-7/2}\) to get \(b_3\) and so on. These \(b_i\) give the right asymptotic
behaviour of \(\hat{G}(k)\) up to a multiplicative constant. To set that multiplicative constant correctly,
we normalise the \((b_1,b_2,\dots, b_p)\) by multiplying a constant so that \(\sum b_i = 2 \). In this way, our near-optimal \(\vec{a}\)
can be used to compute a \(\vec{b}\) which, in turn, gives us a (hopefully near-optimal)
\(G(x)\) satisfying equation~\eqref{eqn_lower_bound_G} and its
Fourier coefficients \(\hat{G}(k)\).

Now, with a good candidate \(\hat{G}(k)\), it remains to compute \(\sum_k
\abs{\hat{G}(k)}^{4/3} \) that is required for
equation~\eqref{eqn_nu_lower_bound}. We proceed in much the same way
we summed \(\hat{F}(k)^4\). We again use Kummer's series transform. So we sum
exactly for small \(n\), and then for large \(n\) we \(\hat{G}(k)\) as per
equation~\eqref{eqn_bess_res}, and then use that to compute the
expansion of its \(4/3\)-power. Again, these expansions exhibit the same
period-4 behaviour, and so we really sum the large-n contributions in four
parts. Putting this all together then gives us the required sum, with rigorous error bounds.

\section{Results and explorations} \label{sec results}
The methods outlined above were programmed in \texttt{c++} using ball-arithmetic facilitated by the
\texttt{flint} library. There are some parameters in these calculations
\begin{itemize}
	\item the number of digits of precision for floating point calculations
	\item \(P\) being the number of coefficients in the ansatz~\eqref{eqn_ansatz_2}
	\item \(N\) being the cutoff between ``small'' and ``large'' \(n\)-values for the Kummers-transform summation~\eqref{eqn_kummer_transform}
	\item \(K\) being the number of terms in the expansions~\eqref{eqn_bess_res}
\end{itemize}
Since all the calculations were carried out multiple times for different
values of the coefficients \(\vec{a}\) we found it helpful to pre-compute and store many parts
of the computations:
\begin{itemize}
	\item the values of the Hurwitz-zeta functions in equation~\eqref{eqn_hurwitz_four} --- used to sum the asymptotic forms.
	\item the values of \(N^k\) for a large range of \(k\)-values --- used to compute and track error terms from truncating expansions.
	\item the values of \(\mybess{p}{\frac{\pi n}{2}}\) for \(0 \leq p \leq P\) and \(1 \leq n \leq N\) --- for small-\(n\) summations.
	\item the asymptotic expansions of \(\mybess{p}{\frac{\pi n}{2}} \) --- used to build expansions such as~\eqref{eqn_bess_res}.
\end{itemize}

\subsection{Main result}

We started with a modest value of \(P=4\), with \(N=1024\) and \(K=64\), 256
digits of floating-point precision, and then set \(\vec{a}=(1,0,0,0),
\lambda=0.5\). The code then uses precomputed expansions of Bessel functions to
assemble the expansion of \(\hat{F}(k)\) and a bound on the truncation error
using the assumption that \(n \geq N\). This is then used to build expansions
of \(\hat{F}(k)^2, \hat{F}(k)^3\) and \(\hat{F}(k)^4\). We then compute the
Lagrangian and the required partial derivatives using Kummers-transform and the
pre-computed Bessel- and zeta-function values, while bounding the error and
using ball-arithmetic to control rigorously the floating-point errors. Much of
this can be parallelised and it typically only took a few seconds on a relatively
modern laptop computer.

The first and second partial-derivatives are used to construct the gradient and Hessian, and then
perform a Newton-Raphson iteration to compute a new value of \(\vec{a}\).
We iterated this process until the norm of the gradient was sufficiently small
--- say, on the order of \(10^{-40}\). For \(P=8\) this was only a few tens of seconds.
At this point, we used \(\vec{a}\) to compute \(\vec{b}\) as outlined above and summed \(\sum
\abs{\hat{G}(k)}^{4/3}\) again using the Kummer-transform. This gives the following upper and lower bounds
on \(\nu_2^2\):
\begin{equation}
	\underline{0.57463\,960}\,65 \leq \nu_2^2 \leq  \underline{0.57463\,960}\,72
\end{equation}
where we have underlined the common digits and have rounded the last digit down
and up respectively. We note that the above bounds are rigorous and any
truncation errors or floating point errors are many orders of magnitude smaller
than the rounding of the last stated digit.

Once we had a near-optimal \(\vec{a}\) for \(P=4\), we increased \(P\) to
\(8\), and then appended \(0\)'s to \(\vec{a}\) to make a new starting point
for Newton-Raphson. We iterated until this had reasonably converged and then
computed upper and lower bounds.
\begin{equation}
	\underline{0.57463\,96071\,51519}\,47 \leq \nu_2^2 \leq  \underline{0.57463\,96071\,51519}\,60
\end{equation}
Again, we underline the common digits and round the last digit down and up respectively. Note that
as \(P\) is increased, the optimal \(\vec{a}\) changes very little; we observe that \(a_j(P)-a_j(P+1)\)
decays exponentially with \(P\).

Then incremented \(P\) again and repeated this process. We found that as \(P\)
became larger it was necessary to increase \(N\) and \(K\) and also the number
of digits of precision. We finished with \(P=101, N=8192, K=128\) and \(384\)
digits of precision. This gives the first 128 digits of \(\nu_2^2\) which is the main result of the paper.

\begin{figure}
	\begin{center}
		\includegraphics[width=0.45\textwidth]{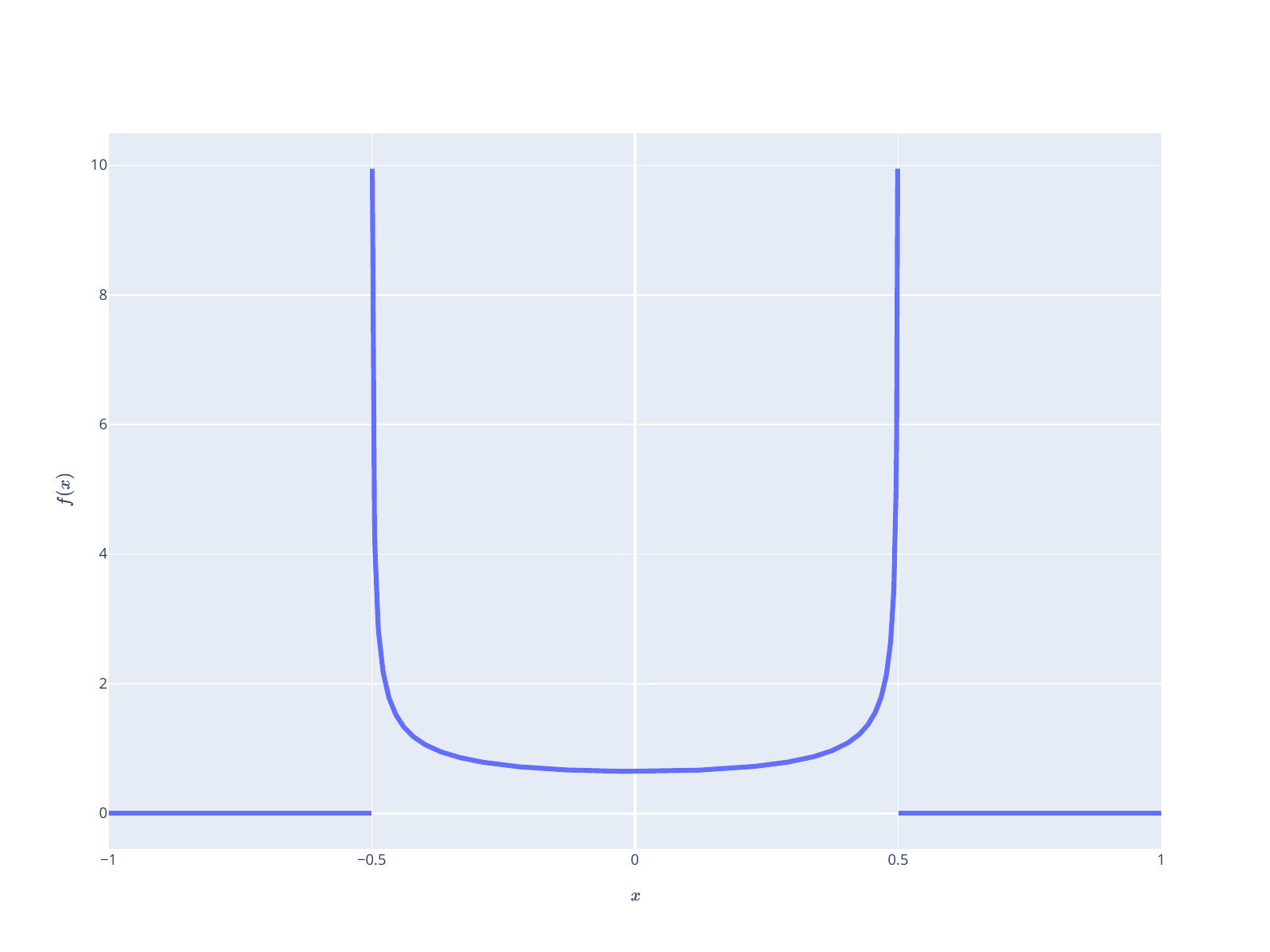}
		\includegraphics[width=0.45\textwidth]{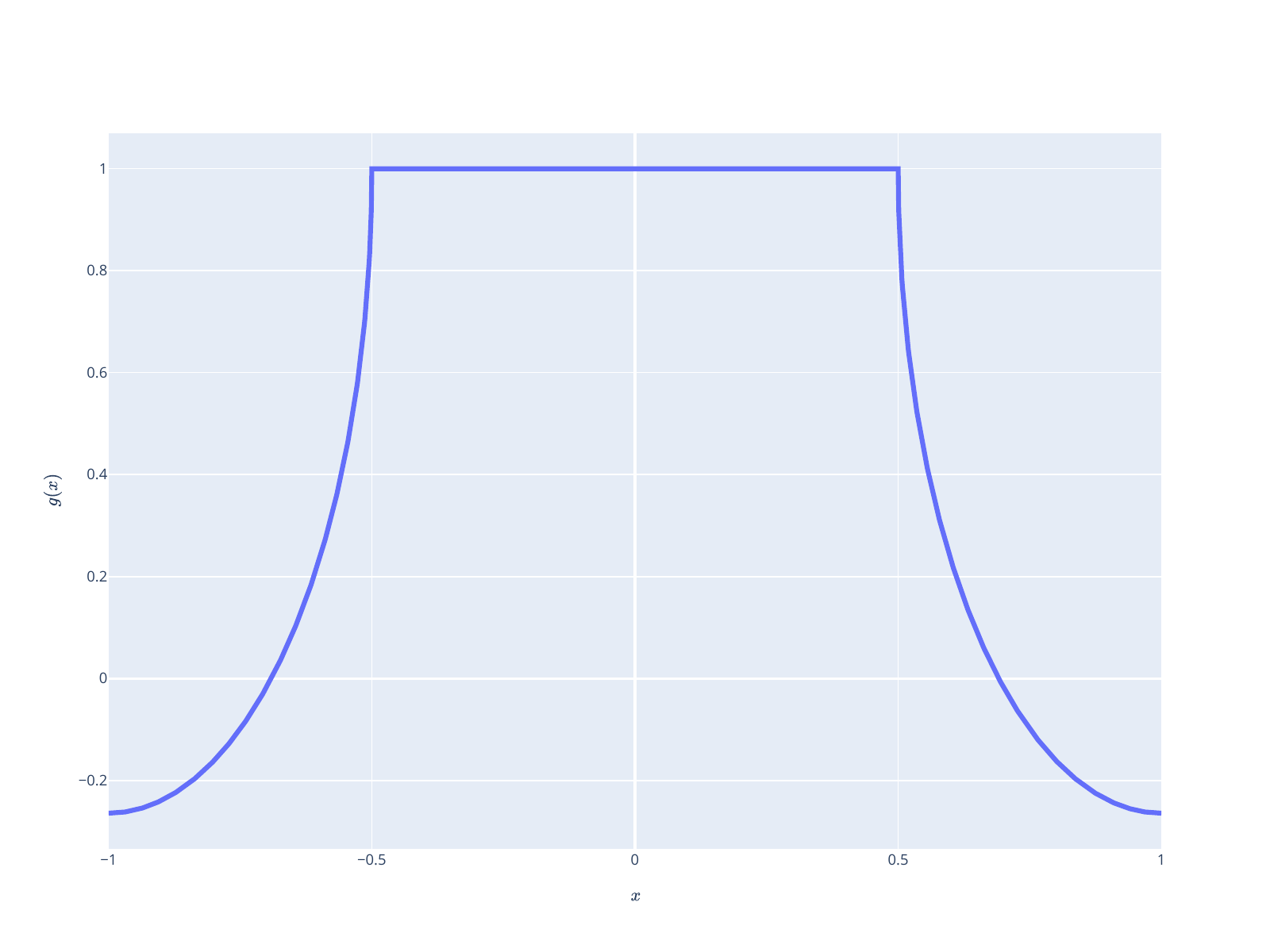}
	\end{center}
	\caption{Plots of the very-nearly optimal \(F(x)\) and \(G(x)\) used to prove Theorem~\ref{thm main}.}
	\label{fig_near_opt_f_and_g}
\end{figure}

\begin{theorem}\label{thm main}
	The quantity \(\nu_2^2 = \inf_{f \in \mathcal{F}} \| f\ast f\|_2^2 \),
	where the infimum is taken over all functions \(f  \in L^1(-1/2,1/2) \) with \(\int f=1\),
	is bounded by
	\[
		c_\ell \leq \nu_2^2 \leq c_u
	\]
	where \( \abs{c_u-c_\ell} \leq 1.2 \times 10^{-129} \), and
	\begin{align*}
		c_\ell & = \underline{
			0.
			57463\,
			96071\,
			51519\,
			59272\,
			72554\,
			27527\,
			05297\,
			14370\,
			26369\,
			37315\,
			66116\,
			30876\, \cdots
		}                                  \\
		       & \qquad \underline{74892\,
			55216\,
			18178\,
			98882\,
			24078\,
			24755\,
			71532\,
			95571\,
			66060\,
			64735\, 74241\, 32638\, 64820\, 673}\, 58
		\\
		c_u    & =
		\underline{
			0.
			57463\,
			96071\,
			51519\,
			59272\,
			72554\,
			27527\,
			05297\,
			14370\,
			26369\,
			37315\,
			66116\,
			30876\, \cdots
		}                                  \\
		       & \qquad \underline{74892\,
			55216\,
			18178\,
			98882\,
			24078\,
			24755\,
			71532\,
			95571\,
			66060\,
			64735\, 74241\, 32638\, 64820\, 673}\, 69
	\end{align*}

	We have underlined first 128 digits that are common to the upper and lower bounds. Consequently, if \(f:[-1/2,1/2] \to \mathbb{R}^+\) is a non-negative function with \( \int f = 1\), then \( \|f \ast f \|_2^2 \geq c_\ell \).
\end{theorem}
The \(\vec{a}\) and \(\vec{b}\) coefficients used to compute the bounds in this result are given in Appendix~\ref{sec cfs}. We also give \texttt{python} code in Appendix~\ref{sec code} that demonstrates how to use
use (the first few of) these coefficients to compute upper and lower bounds --- albeit non-rigorously. We note that while the \(\vec{a}\)-coefficients are steadily decreasing in magnitude, the \(\vec{b}\)-coefficients stabilise in magnitude at about \(b_70\). This suggests that fewer \(b\)-coefficients might be necessary to compute the lower bound to the given accuracy; we have not explored this possibility. In Figure~\ref{fig_near_opt_f_and_g} we plot the corresponding (very near) optimal \(F(x) \) and \(G(x)\) that are used to prove Theorem~\ref{thm main}.

\begin{figure}
	\begin{center}
		\includegraphics[width=0.66\textwidth]{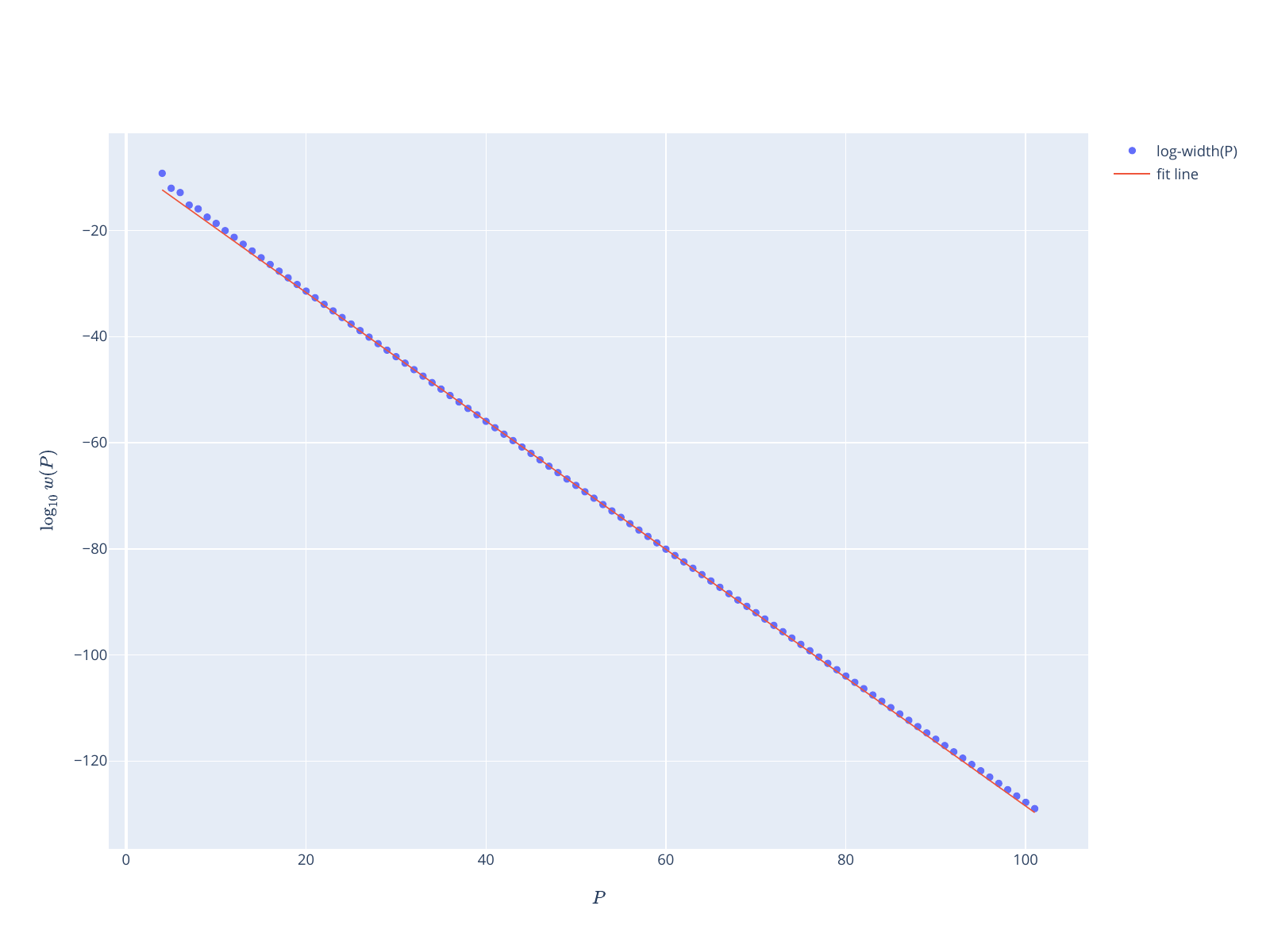}
	\end{center}
	\caption{The decimal-logarithm of the width of the bound against \(P\). We also show a simple linear fit of the data
	which shows that \(\log_{10}(w(P)) \approx -7.5 -1.2P \). }\label{fig_bound_width}
\end{figure}

We note that there is no real technical impediment to extending these bounds to more digits, other than computer time and patience; the memory requirements were negligible.
In Figure~\ref{fig_bound_width} we have plotted the width of the bound (i.e. upper minus lower) against \(P\), the number of coefficients in the ansatz.
We see that the logarithm of the width of the bound is quite linear against \(P\). A simple rough linear fit gives \(\log_{10} w(P) \approx -7.6-1.2P\)
and so we extrapolate that we would need about \(P\approx 160\) to fix the first 200 digits and \(P \approx 830\) to get the first
1000 digits.

\subsection{Explorations}
Part of the motivation of pushing the above calculations to this level of
precision was the hope that we could find some underlying closed form
expression for \(\nu_2\) or the corresponding minimising function. Unfortunately
despite the precision with which we now know \(\nu_2^2\) we have not been able
to guess a simple closed form using RIES~\cite{ries}, the Inverse Symbolic
Calculator or Plouffe's inverter~\cite{stoutemyer2021}.

We have also examined the \(\vec{a}\) coefficients, hoping to guess some closed
forms. We do find that the coefficients converge quickly as \(P\) increases,
and we also observe that
\begin{equation}
	\frac{a_{n}}{a_{n+1}} \sim -8 \cdot (1 + n^{-1}).
\end{equation}
This simple form seems to hold except for \(n\)-values close to \(P\). Similar behaviour was observed for smaller values of \(P\) and
we expect that if \(P\) were increased further, that the above functional form would continue to hold.
This, in turn suggests that \(a_n \sim A (-8)^{-n} n^{-1} (1 + \dots ) \). Some rough fitting yielded
\begin{equation}
	n (-8)^n a_n \approx 1.45 + 0.54 n^{-1},
\end{equation}
but could not find any particularly convincing values of the above constants, but we think that they give an indication of the magnitude of those terms. That
said, the signal of \( (-8)^{-n} n^{-1}\) in the coefficients is extremely robust.

\begin{figure}
	\begin{center}
		\includegraphics[width=0.66\textwidth]{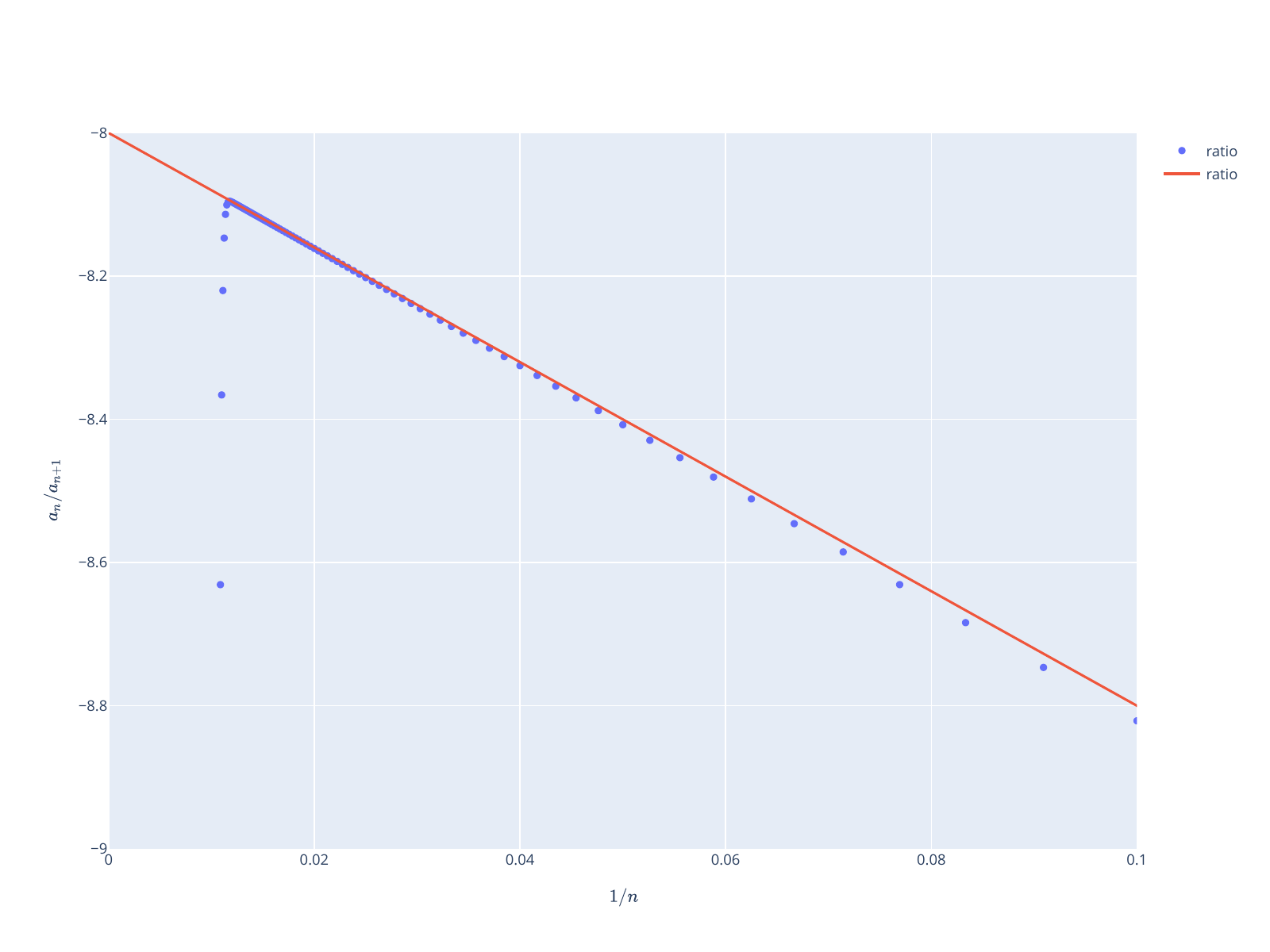}
	\end{center}
	\caption{A plot of the ratios \(a_{n}/a_{n+1}\) against \(1/n\) obtained with \(P=101\). We note that the data is extremely well fitted by the
	line \(y=-8\cdot (1+1/n)\) until \(n\) is close to the final few coefficients. }\label{fig_cfs_ratio}
\end{figure}

We note that we see similar scaling of the \(\vec{b}\)-coefficients used to compute the lower bound:
\begin{equation}
	\frac{b_{n}}{b_{n+1}} \sim -8 \cdot (1 + 2 n^{-1}).
\end{equation}
Again, while the signal of \((-8)^{-n} n^{-2}\) is very strong and robust, we did not arrive at any convincing
fit of more detailed asymptotics. Very roughly we observe
\begin{equation}
	n^2 (-8)^n b_n \approx 4.5 - 2.3 n^{-1},
\end{equation}
but this should be taken as an indication of the magnitude of those coefficients, not a reliable fit.

We are not entirely sure what to make of these asymptotic signals, other than
that it suggests there is more structure in this problem waiting to be
uncovered.

\section*{Acknowledgements}
The author thanks Michael Bennett, Tony Guttmann, Kevin O'Bryant, Jozsef
Solymosi and Ethan White for their helpful discussions. I would also like
to give credit to the developers of \texttt{flint} and \texttt{mpmath} for
their excellent software; it was key both for numerical explorations and rigorous calculations.
The author also acknowledges funding from NSERC via the Discovery Project program.

\printbibliography

@book{abramowitz1964,
  title = {Handbook of mathematical functions with formulas, graphs, and
           mathematical tables},
  author = {Abramowitz, Milton and Stegun, Irene A},
  volume = {55},
  year = {1964},
  publisher = {US Government printing office},
}

@article{cilleruelo2000,
  title = {{\(B_h[g]\)} sequences},
  author = {Cilleruelo, Javier and Jim{\'e}nez-Urroz, Jorge},
  journal = {Mathematika},
  volume = {47},
  number = {1-2},
  pages = {109--115},
  year = {2000},
  publisher = {London Mathematical Society},
}

@article{cilleruelo2008,
  title = {{\(B_2[g]\)} sets and a conjecture of Schinzel and Schmidt},
  author = {Cilleruelo, Javier and Vinuesa, Carlos},
  journal = {Combinatorics, Probability and Computing},
  volume = {17},
  number = {6},
  pages = {741--747},
  year = {2008},
  publisher = {Cambridge University Press},
}

@article{cilleruelo2010,
  title = {Generalized sidon sets},
  author = {Cilleruelo, Javier and Ruzsa, Imre and Vinuesa, Carlos},
  journal = {Advances in Mathematics},
  volume = {225},
  number = {5},
  pages = {2786--2807},
  year = {2010},
  publisher = {Elsevier},
}

@article{cloninger2017,
  title = {On suprema of autoconvolutions with an application to Sidon sets},
  author = {Cloninger, Alexander and Steinerberger, Stefan},
  journal = {Proceedings of the American Mathematical Society},
  volume = {145},
  number = {8},
  pages = {3191--3200},
  year = {2017},
}

@article{erdos1941,
  author = {Erd{\H{o}}s, Paul and Tur{\'a}n, Paul},
  title = {On a problem of Sidon in additive number theory, and on some related
           problems.},
  joural = {J. London Math. Soc},
  year = {1941},
  volume = 16,
  issue = 4,
  pages = {212--215},
}

@book{fletcher2013,
  title = {Practical methods of optimization},
  author = {Fletcher, Roger},
  year = {2013},
  publisher = {John Wiley \& Sons},
}

@manual{flint,
  key = {{FLINT}},
  author = {The {FLINT} team},
  title = {{FLINT}: {F}ast {L}ibrary for {N}umber {T}heory},
  year = {2025},
  note = {Version 3.4.0, \\ \url{https://flintlib.org}},
}

@article{georgiev2025,
  title = {Mathematical exploration and discovery at scale},
  author = {Georgiev, Bogdan and G{\'o}mez-Serrano, Javier and Tao, Terence and
            Wagner, Adam Zsolt},
  journal = {arXiv: 2511.02864},
  year = {2025},
}

@article{green2001,
  title = {The number of squares and $ B_h [g] $ sets},
  author = {Green, Ben},
  journal = {Acta Arithmetica},
  volume = {100},
  number = {4},
  pages = {365--390},
  year = {2001},
}

@misc{green100,
  title = {100 open problems},
  author = {Green, Ben},
  note = {manuscript, available on request to Professor Green},
}

@article{johansson2017,
  author = {Fredrik Johansson},
  journal = {{IEEE} Transactions on Computers},
  title = {{A}rb: Efficient arbitrary precision midpoint radius interval
           arithmetic},
  year = {2017},
  volume = {66},
  number = {8},
  pages = {1281-1292},
}

@article{johnston2022,
  title = {Upper and lower bounds on the size of {\(B_k[g]\)} sets},
  author = {Johnston, Griffin and Tait, Michael and Timmons, Craig},
  journal = {Australiasian Journal of Combinatorics},
  volume = {83},
  number = {1},
  pages = {129--140},
  year = {2022},
}

@article{lindstrom2000,
  title = { {\(B_h[g]\)}-sequences from {\(B_h\)}-sequences},
  author = {Lindstr{\"o}m, Bernt},
  journal = {Proceedings of the American Mathematical Society},
  pages = {657--659},
  year = {2000},
  publisher = {JSTOR},
}

@article{martin2007,
  title = {The symmetric subset problem in continuous Ramsey theory},
  author = {Martin, Greg and O'Bryant, Kevin},
  journal = {Experimental Mathematics},
  volume = {16},
  number = {2},
  pages = {145--165},
  year = {2007},
  publisher = {Taylor \& Francis},
}

@article{martin2009,
  title = {The supremum of autoconvolutions, with applications to additive
           number theory},
  author = {Martin, Greg and O’Bryant, Kevin},
  journal = {Illinois Journal of Mathematics},
  volume = {53},
  number = {1},
  pages = {219--235},
  year = {2009},
  publisher = {Duke University Press},
}

@article{obraynt2004,
  title = {A complete annotated bibliography of work related to Sidon sequences},
  author = {O'Bryant, Kevin},
  journal = {Electronic Journal of Combinatorics},
  volume = {Dynamic surveys 11},
  pages = {1--39},
  year = {2004},
}

@manual{mpmath,
  key = {mpmath},
  author = {The mpmath development team},
  title = {mpmath: a {P}ython library for arbitrary precision floating point
           arithmetic (version 1.3.0)},
  note = {\url{http://mpmath.org/}},
  year = {2023},
}

@misc{olver2019,
  title = {Digital library of mathematical functions},
  author = {Olver, Frank W~J and Olde Daalhuis, Adri B and Lozier, Daniel W and
            Schneider, Barry and Boisvert, Ronald F and Clark, Charles W and
            Miller, Bruce R and Saunders, Bonita V and Cohl, Howard S and McClain
            , Marjorie A},
  year = {2019},
  publisher = {National Institute of Standards and Technology (NIST)},
  note = {\url{https://dlmf.nist.gov/}},
}

@book{olver1997,
  title = {Asymptotics and special functions},
  author = {Olver, Frank},
  year = {1997},
  publisher = {AK Peters/CRC Press},
}

@misc{ries,
  author = {Robert Munafo},
  title = {Ries: Rilybot Inverse Equation Solver},
  year = {2022},
  note = {\url{https://www.mrob.com/pub/ries/index.html}},
}

@article{schinzel2002,
  title = {Comparison of {$L^{1}$}- and {$L^{\infty}$}-norms of squares of
           polynomials},
  author = {Schinzel, A and Schmidt, WM},
  journal = {Acta Arithmetica},
  volume = {104},
  pages = {283--296},
  year = {2002},
  publisher = {Instytut Matematyczny Polskiej Akademii Nauk},
}

@book{sidi2003,
  title = {Practical extrapolation methods: Theory and applications},
  author = {Sidi, Avram},
  volume = {10},
  year = {2003},
  publisher = {Cambridge university press},
}

@article{sidon1932,
  title = {Ein Satz {\"u}ber trigonometrische Polynome und seine Anwendung in
           der Theorie der Fourier-Reihen},
  author = {Sidon, Simon},
  journal = {Mathematische Annalen},
  volume = {106},
  number = {1},
  pages = {536--539},
  year = {1932},
  publisher = {Springer},
}

@article{stoutemyer2021,
  title = {How to hunt wild constants},
  author = {Stoutemyer, David R},
  journal = {arXiv: 2103.16720},
  year = {2021},
}

@article{white2022,
  title = {An almost-tight $ L^2$ autoconvolution inequality},
  author = {White, Ethan Patrick},
  journal = {arXiv: 2210.16437},
  year = {2022},
}

@article{white2024,
  title = {An optimal autoconvolution inequality},
  author = {White, Ethan Patrick},
  journal = {Canadian Mathematical Bulletin},
  volume = {67},
  number = {1},
  pages = {108--121},
  year = {2024},
  publisher = {Canadian Mathematical Society},
}

@misc{white_data,
  author = "White, Ethan Patrick",
  date = "2024",
  howpublished = "personal communication",
}

@article{yuksekgonul2026,
  title = {Learning to Discover at Test Time},
  author = {Yuksekgonul, Mert and Koceja, Daniel and Li, Xinhao and Bianchi,
            Federico and McCaleb, Jed and Wang, Xiaolong and Kautz, Jan and Choi,
            Yejin and Zou, James and Guestrin, Carlos and others},
  journal = {arXiv: 2601.16175},
  year = {2026},
}

\appendix
\begin{landscape}
	\section{Ansatz coefficients} \label{sec cfs}
	We give the \(\vec{a}\) and \(\vec{b}\) coefficients that are used to compute the bounds in Theorem~\ref{thm main}. We have truncated these at 128 digits.
		{\scriptsize \ttfamily
			\begin{longtable}{|l|r|}
				\hline
				\(p\) & \(a_p\)                                                                                                                             \\[0.5ex]
				\hline \endhead
				\hline
				\(p\) & \(a_p\)                                                                                                                             \\[0.5ex]
				\hline \endfoot
				0     & 0.98602065571593947981982597405697290479774950169551574408171685338064128443650486639695959919838106504681745070718453265042295595  \\
				1     & 0.01489770489081501391236554516465932337360037886185262531409922769587507544257885553721071292939933006377392624509980008292575019  \\
				2     & -0.00100205147026116189707339482723743634634855852968990358063334091805583926892625550551551017653327044507577012139566932195157950 \\
				3     & 0.00009204407950297445359544529687886064943150128923593264098550129302924478408374955135440082953010229598917584177918942923344904  \\
				4     & -0.00000922274904509772914881249306215878382678034496365492938748571679130161109502434346988480188665105871343249815026801667610920 \\
				5     & 0.00000096248293777695239529770940406980385573510172484739923996188253946427333183826276677870868029861653768632831444028178007925  \\
				6     & -0.00000010308254667112377770352811050554776487355578124053305283669721149458702219139700107199486390915352511624730946494476325852 \\
				7     & 0.00000001125442999323010976782605876372232519419601097844133358508048774936750593188528772457877750953930068098950646291457515379  \\
				8     & -0.00000000124752688834370157801836354986319635941999009201009510429123187129057108187123588457671365713289942497871930592821722596 \\
				9     & 0.00000000013999835706501306249491275040186120544384912063038540963472409378032308691716797472525682861221215931077648801243110953  \\
				10    & -0.00000000001587108251776074784839022898212750147053313531428602050340190360631660966934550589405568798555200352106118884222032538 \\
				11    & 0.00000000000181457936724973745938016642332925399974993453212544266259390429885622289685497037488525637549662968871641547375478972  \\
				12    & -0.00000000000020895936361219667407322629723937158854300531603277784230504175662574715902488555569213215274038232249436834427175533 \\
				13    & 0.00000000000002421093838520663620792232232117921974355045462840913540269329904943007360029724462855163906982593109397378696703991  \\
				14    & -0.00000000000000282008054061566739367389444663724772159543038013815771878818399242827152032103549753349436868166270269801819363472 \\
				15    & 0.00000000000000033000131420654058434725678872629262934997700068966298963247713780231593921874957824445959306112754673467111387467  \\
				16    & -0.00000000000000003877312975807208249544469859456112800747365135241934148282443154017384602330494086694938790432875367855607117426 \\
				17    & 0.00000000000000000457198691425878853747411301567813802192928752567789475706546885071732137722406542546999659597636047872144377889  \\
				18    & -0.00000000000000000054083956779454257619021646503718481432293544952313528951804923049437792009155201232560572092281770582103834112 \\
				19    & 0.00000000000000000006416197485660862347940887486845107212042017100285527469759052073014019639801219350633648958205415722289264059  \\
				20    & -0.00000000000000000000763150526231448395786846231154147604545362252211343391495743473317161110587576315905247644179588270741479441 \\
				21    & 0.00000000000000000000090983047251504819787197629681843879759005992127979856128680460179042780177735369732448841924593108178994189  \\
				22    & -0.00000000000000000000010870191728598186414407822486181171595210052896685470598994002348238365199433796607771128381837768217912410 \\
				23    & 0.00000000000000000000001301250076075800810441815260661785624451441008797432775406046786444831381316526909774699990897544927610652  \\
				24    & -0.00000000000000000000000156049130477462767383517075901464180143779239226022517349114356840077165476859649170862277620551175965633 \\
				25    & 0.00000000000000000000000018744653061723895081638820842105919591953375602299779714965620404571655662221458572173057040335037646517  \\
				26    & -0.00000000000000000000000002255040723271186408025465369462058839681249458584143153990360285002830796325441774423267289774109646182 \\
				27    & 0.00000000000000000000000000271671294242622470411133487840687989946022016102614595674805919116951699405643831033233150955970035971  \\
				28    & -0.00000000000000000000000000032771935591411601308066385942808429522025085672859175030664123613989878132858385226641230329614128915 \\
				29    & 0.00000000000000000000000000003958133665045999892178660686839121277569959248873313737614558281165992351168678903542940526815320742  \\
				30    & -0.00000000000000000000000000000478601159937482790532105231810477017342000309129188325857890873035056021031286321243340010501132973 \\
				31    & 0.00000000000000000000000000000057932228794767412709441362326582287633563383939790182510297892574313694959864399612127172753341730  \\
				32    & -0.00000000000000000000000000000007019418737404626005609577331569335732864594672763388689711123571815852534728551550610361359368438 \\
				33    & 0.00000000000000000000000000000000851315042656325612270758598959076419037368451980922296455814784658427453908722736741057688772015  \\
				34    & -0.00000000000000000000000000000000103338896966892030571003062238739882733489585104515513691361635797928809560274284474059058823274 \\
				35    & 0.00000000000000000000000000000000012554512150951189429585281121201242465890427601060360317610528088904704412980452539092813718535  \\
				36    & -0.00000000000000000000000000000000001526435031665054994534643884594676617003584343451307931585740959599905812727909767798748019287 \\
				37    & 0.00000000000000000000000000000000000185729480712076063363385226660885583492923010917134772323815446731339740028162825812204455558  \\
				38    & -0.00000000000000000000000000000000000022614678667523558434996310714672335220277935561460879914156266652285326045728268389137715257 \\
				39    & 0.00000000000000000000000000000000000002755442635583815710739176959606588537226473609804094562343881209544923164259406997861163794  \\
				40    & -0.00000000000000000000000000000000000000335945781811372694074856956815985583602687461844042941958670005067893881321666238027815746 \\
				41    & 0.00000000000000000000000000000000000000040983630019698508285790704521302174787362386873246780547455289544119912211694645825592274  \\
				42    & -0.00000000000000000000000000000000000000005002678186810020583539468276657847610853800940511256232735582198524246062313541324029721 \\
				43    & 0.00000000000000000000000000000000000000000610989849801212933800736986084777349075949075376896244314141085524033827057338656651215  \\
				44    & -0.00000000000000000000000000000000000000000074661011576241763781149696564640636011206071363289428811286822317218436420886111916034 \\
				45    & 0.00000000000000000000000000000000000000000009127925136644309830243688356670379927578270041714958780758124317696714360299671667879  \\
				46    & -0.00000000000000000000000000000000000000000001116501040130585035565237462382398547342129876462018316308073098317657749869089398504 \\
				47    & 0.00000000000000000000000000000000000000000000136630040314889418946084111095472357974241937838038854566595124021717476668282637379  \\
				48    & -0.00000000000000000000000000000000000000000000016727265646568911299231591021763287414204923035851678693220709379767167768597450686 \\
				49    & 0.00000000000000000000000000000000000000000000002048743565445389116475451328962090191761415085063032706097210712842041946377730919  \\
				50    & -0.00000000000000000000000000000000000000000000000251030682428190474968916660237616098011845620963629883553284804565875006775962962 \\
				51    & 0.00000000000000000000000000000000000000000000000030770578488894628208517329488323662091571197928212443912914206930105027858417610  \\
				52    & -0.00000000000000000000000000000000000000000000000003773181161746214009070449857978871987483620956858431708996825918865849511269481 \\
				53    & 0.00000000000000000000000000000000000000000000000000462846149360105202002484177171659644203262883969784426174575427213517702165059  \\
				54    & -0.00000000000000000000000000000000000000000000000000056795882600461820374531327015489786234964653194074467627260877757274923746785 \\
				55    & 0.00000000000000000000000000000000000000000000000000006971765417249870578431120081035902332169902580272201971670097651986761089733  \\
				56    & -0.00000000000000000000000000000000000000000000000000000856069849426285765919176101721565535014743374907968143395936375000981364096 \\
				57    & 0.00000000000000000000000000000000000000000000000000000105150475428056475261872383226343216452461168333790566472181993704523484489  \\
				58    & -0.00000000000000000000000000000000000000000000000000000012919455700619248916665905455588455130410838120521596875829801567317032934 \\
				59    & 0.00000000000000000000000000000000000000000000000000000001587828898393624011661993251009053574233906022872417385501389485525302844  \\
				60    & -0.00000000000000000000000000000000000000000000000000000000195202549878205616820430361201020206702181094495096006654193338334890111 \\
				61    & 0.00000000000000000000000000000000000000000000000000000000024004107840376710474169451740849130355531835352818147826397737584550563  \\
				62    & -0.00000000000000000000000000000000000000000000000000000000002952569583647896352995075016457509102074807989029557375349591707840394 \\
				63    & 0.00000000000000000000000000000000000000000000000000000000000363266689081628479448378028542451918607556164994539452242570579105826  \\
				64    & -0.00000000000000000000000000000000000000000000000000000000000044705238245710279823577791099792192493580614814526490415429879814107 \\
				65    & 0.00000000000000000000000000000000000000000000000000000000000005502947340108063904373001898010057993913926859547534443917891125410  \\
				66    & -0.00000000000000000000000000000000000000000000000000000000000000677537373390466410062613440224799220415853549235538350548111060287 \\
				67    & 0.00000000000000000000000000000000000000000000000000000000000000083439002422714847661370691263451003815932784444023177656577044147  \\
				68    & -0.00000000000000000000000000000000000000000000000000000000000000010277796675991542873549909134179999472113333055760273964752408880 \\
				69    & 0.00000000000000000000000000000000000000000000000000000000000000001266261087869095894645063602792465850594025158906214580661149616  \\
				70    & -0.00000000000000000000000000000000000000000000000000000000000000000156040088452232213271828091877736247139627119146804989306153944 \\
				71    & 0.00000000000000000000000000000000000000000000000000000000000000000019232523831051335348807208205304868932362345547978385313398376  \\
				72    & -0.00000000000000000000000000000000000000000000000000000000000000000002370943026407117817636477802257088400739853564199485451725501 \\
				73    & 0.00000000000000000000000000000000000000000000000000000000000000000000292340096447487567327508866641171738267669130416855354393651  \\
				74    & -0.00000000000000000000000000000000000000000000000000000000000000000000036052539253932852363452218059569604866043306342254080023117 \\
				75    & 0.00000000000000000000000000000000000000000000000000000000000000000000004446941331165304635929019689784462789499993745600082553745  \\
				76    & -0.00000000000000000000000000000000000000000000000000000000000000000000000548609038770849939160083300295769952024471788809309455791 \\
				77    & 0.00000000000000000000000000000000000000000000000000000000000000000000000067692189743857336721825798685418520282524305563802493393  \\
				78    & -0.00000000000000000000000000000000000000000000000000000000000000000000000008353843138161365245490941488596865498932467542025236714 \\
				79    & 0.00000000000000000000000000000000000000000000000000000000000000000000000001031108551154953722565775350131044229825405436427577762  \\
				80    & -0.00000000000000000000000000000000000000000000000000000000000000000000000000127289036532278537174304773357456685108072477703753884 \\
				81    & 0.00000000000000000000000000000000000000000000000000000000000000000000000000015716084067912821226818145169287808062622907570359639  \\
				82    & -0.00000000000000000000000000000000000000000000000000000000000000000000000000001940717446638463468610019274379634503406862702264311 \\
				83    & 0.00000000000000000000000000000000000000000000000000000000000000000000000000000239685326985640090722064041603253236123708795368800  \\
				84    & -0.00000000000000000000000000000000000000000000000000000000000000000000000000000029605571640403838995092582416190811408531961256608 \\
				85    & 0.00000000000000000000000000000000000000000000000000000000000000000000000000000003657078619107130271227666881993110356454449313318  \\
				86    & -0.00000000000000000000000000000000000000000000000000000000000000000000000000000000451704063474421387286818571796425847890279116776 \\
				87    & 0.00000000000000000000000000000000000000000000000000000000000000000000000000000000055762414595158347063192888926077011963403074360  \\
				88    & -0.00000000000000000000000000000000000000000000000000000000000000000000000000000000006872615798301534950869865375959585368851772949 \\
				89    & 0.00000000000000000000000000000000000000000000000000000000000000000000000000000000000843593693224077402456720589251805983549672198  \\
				90    & -0.00000000000000000000000000000000000000000000000000000000000000000000000000000000000102625981677427532379241020055915482302818245 \\
				91    & 0.00000000000000000000000000000000000000000000000000000000000000000000000000000000000012267176568283274837672614763384918773986133  \\
				92    & -0.00000000000000000000000000000000000000000000000000000000000000000000000000000000000001421319507618214156346648325152478178549807 \\
				93    & 0.00000000000000000000000000000000000000000000000000000000000000000000000000000000000000156580949030229658683142149360551107251072  \\
				94    & -0.00000000000000000000000000000000000000000000000000000000000000000000000000000000000000015995318062396177094521184759359302325079 \\
				95    & 0.00000000000000000000000000000000000000000000000000000000000000000000000000000000000000001468661860210306721310183054595937895505  \\
				96    & -0.00000000000000000000000000000000000000000000000000000000000000000000000000000000000000000116614925472654571330498298650548807495 \\
				97    & 0.00000000000000000000000000000000000000000000000000000000000000000000000000000000000000000007614777707935875830049462101348949930  \\
				98    & -0.00000000000000000000000000000000000000000000000000000000000000000000000000000000000000000000380292841621387853386840208666017802 \\
				99    & 0.00000000000000000000000000000000000000000000000000000000000000000000000000000000000000000000012827656144638619781338571732916502  \\
				100   & -0.00000000000000000000000000000000000000000000000000000000000000000000000000000000000000000000000218049257000709017436253012230237
			\end{longtable}
		}

		{\scriptsize \ttfamily
			\begin{longtable}{|l|r|}
				\hline
				\(p\) & \(b_p\)                                                                                                                             \\[0.5ex]
				\hline \endhead
				\hline
				\(p\) & \(b_p\)                                                                                                                             \\[0.5ex]
				\hline \endfoot
				1     & 1.02206573731609905553110432296465902673166228639065467690912389968916637683205223973582386543900423140398430546973159699475880968  \\
				2     & -0.02316346059073084765116552642720018104545001014596264722478965457659723311382886327186557407279722375919299170813555292706053273 \\
				3     & 0.00117055009738371923883586925558687165622850288142052568396302488019747209991589548600529080807443850935071498170077042683798424  \\
				4     & -0.00007845159525437915038502792165832119901038961215248844398808658943694139861261149529979387302170261581669509813473783075755960 \\
				5     & 0.00000610178237044628533230298822268990953790262896190017750194403093348467248340832520887861010195294760869700709601084842958146  \\
				6     & -0.00000052009749369734457624196511856182756842364965190922715873326904011722795939095756133213413459772421016866628709339695479837 \\
				7     & 0.00000004716048719323026016962146240115292301730419351360452195728031621509007487530119784377060607512711875733268720745381058157  \\
				8     & -0.00000000447126480926096747352520489301311549059663732590268047760546569224386977254240244838688643305230022562630847483191266594 \\
				9     & 0.00000000043843194080269559182644754603387467792414053175117043821909145375871527823423317212393269059722078626700962663261736071  \\
				10    & -0.00000000004413778800416697987488111621886316972971777963650253306012024129999288565136520457054049956121834345143245252834524094 \\
				11    & 0.00000000000453850435772300502954620138881926113395849977050079238233303792339074334569809196529392701461789312735380240119972256  \\
				12    & -0.00000000000047486278974912034286589067491552878774479637193124576808690075176362156330115195535199055003415188530210905162180252 \\
				13    & 0.00000000000005041269366573337264284738549684527612711308685656056720785968677824796325858875903372436174310168374512545650535332  \\
				14    & -0.00000000000000541841618794557231371070940216068718266116912845079186467540276563822875049011184703554378702570907294736853128573 \\
				15    & 0.00000000000000058858832900333091258330415756181530252816126430130927258440399742988060678786813627227959312721665247343200197643  \\
				16    & -0.00000000000000006452865046019881690193229702423234394266581730429386093120750533882537148265525101956496103051609251347741235428 \\
				17    & 0.00000000000000000713184602745285983051433421344021465203851270055667989477399624009688486431914661536929826202761724811871840587  \\
				18    & -0.00000000000000000079387568585334837070391092205296679510736843569182769646154851948665023095367011771242662073488945666936998953 \\
				19    & 0.00000000000000000008893321548244486214039081944855655945918195093438781586307498990499528098749606669065758926471710467540047864  \\
				20    & -0.00000000000000000001001957549700270083508893225362985836632615526906939557708140345148409323539993515199249935730116780451406997 \\
				21    & 0.00000000000000000000113465531831668981697562638356703154556024045986347218789808664711742215819338248894742120026825744832313722  \\
				22    & -0.00000000000000000000012909158102361956404235943698783179138024781235454548651727492369987543101334262582915647238742305585314009 \\
				23    & 0.00000000000000000000001474923829230891569573392420268330241574156683436172535947440264730427883175637259680058394418040445465534  \\
				24    & -0.00000000000000000000000169169006107069139629791161738387615344718638919176357414660585878085398954275870491626263567271342893412 \\
				25    & 0.00000000000000000000000019472106989649246522258324010344258833344100105523240948613184828934792986360979817277966652468538217674  \\
				26    & -0.00000000000000000000000002248660696695523188014261533893506417678074602881156248759550802280384621121670987985737153626023553731 \\
				27    & 0.00000000000000000000000000260462735465193754939342701335913646657050196128585470097160099795012169041977526170349037289687768861  \\
				28    & -0.00000000000000000000000000030253924485017718865960871236635856160042902739053621928855674487536279588648581836346115874581093123 \\
				29    & 0.00000000000000000000000000003523272909337653237779542395823859099308870029672548617048433030539547793520222875984603428230495978  \\
				30    & -0.00000000000000000000000000000411303223137921622816749678238895201096694999113287940394757075059166303156268110629856365016520934 \\
				31    & 0.00000000000000000000000000000048123774047691143253802517132972575397282638184456385665934728886394624805178954746991056266709423  \\
				32    & -0.00000000000000000000000000000005642559349469213203367668939112137942851003184997428253527781183666879008568436602935012970315863 \\
				33    & 0.00000000000000000000000000000000662909868573311871789588471291109933754192484601934454868859162683920749247481661582702308554391  \\
				34    & -0.00000000000000000000000000000000078026633675040835701852286162686108717083549428914666313976341929943958384902663131936751812823 \\
				35    & 0.00000000000000000000000000000000009200130410329765264938462253338931973079948091781780377169263455277592349130986816224392341534  \\
				36    & -0.00000000000000000000000000000000001086587088795343412248534254116285808847177753199971566014988266311228544350287691735277320142 \\
				37    & 0.00000000000000000000000000000000000128533058235929497049089952492775974975139938117167349848775072258580679482823863801885806732  \\
				38    & -0.00000000000000000000000000000000000015226790894999180428671756131228272201703708773282754084747638362465157337376704629005005967 \\
				39    & 0.00000000000000000000000000000000000001806390041204563183044450205988022085553421135844580500614529801585302142501674671979207257  \\
				40    & -0.00000000000000000000000000000000000000214581963070086979607000993486744927329290148261618641156268026953663199927729947014483562 \\
				41    & 0.00000000000000000000000000000000000000025522586505921162949136824540185120502549546450899750990168464760717663289547368606697893  \\
				42    & -0.00000000000000000000000000000000000000003039339976397793471210448164934553750192373732833177466011430972211563321386531014706012 \\
				43    & 0.00000000000000000000000000000000000000000362353317531057898850627720443445721610579026548071159561650923261004232738730627611801  \\
				44    & -0.00000000000000000000000000000000000000000043247450991141856737163055563821568362897913470041436822399776951815884914383355968842 \\
				45    & 0.00000000000000000000000000000000000000000005167049445197177381919496760235514626182986883177659857585900770340527432742853727039  \\
				46    & -0.00000000000000000000000000000000000000000000617957330489485490032415101169389655647864765492672800562003867200892279845821312894 \\
				47    & 0.00000000000000000000000000000000000000000000073975748664476103535099888790133409563430535398366541775483950735578690106644027440  \\
				48    & -0.00000000000000000000000000000000000000000000008863753948463911575245954680926253761937997338041899517984896295978772884259738918 \\
				49    & 0.00000000000000000000000000000000000000000000001062984621224640716318940995956490914693495781109767212517284433790115286550583921  \\
				50    & -0.00000000000000000000000000000000000000000000000127585631090748460335848430470852384615778486097094289869491803958259876008973712 \\
				51    & 0.00000000000000000000000000000000000000000000000015325953316364338251976393341603765524212295949584887398909600348773564440122665  \\
				52    & -0.00000000000000000000000000000000000000000000000001842427816163593516505308530593604446333112679410957856954037338130729231310612 \\
				53    & 0.00000000000000000000000000000000000000000000000000221655140079606131603702873412117410068161214605032572141466107747405482737792  \\
				54    & -0.00000000000000000000000000000000000000000000000000026685619615272784141429265485048830685691060109087240779736517357315760098456 \\
				55    & 0.00000000000000000000000000000000000000000000000000003214973712363550244199448444032422240545996397390073265661769536183709113685  \\
				56    & -0.00000000000000000000000000000000000000000000000000000387585327835244441098543216619581741739710274249376702318687193208162913564 \\
				57    & 0.00000000000000000000000000000000000000000000000000000046755929424754938169215207020422232023101579971140376456833241676591775007  \\
				58    & -0.00000000000000000000000000000000000000000000000000000005643853879849716055792610516651871494725212852846988413720937526287923825 \\
				59    & 0.00000000000000000000000000000000000000000000000000000000681671706622268944664004713816863763326391262449156861756806499728842564  \\
				60    & -0.00000000000000000000000000000000000000000000000000000000082380879646684289067620231937844894224228829618260088509451981960933605 \\
				61    & 0.00000000000000000000000000000000000000000000000000000000009961411403265835887411781837110664340902329075041693576559242917464528  \\
				62    & -0.00000000000000000000000000000000000000000000000000000000001205176568507052529041691720167828605605177530275339794336283572721083 \\
				63    & 0.00000000000000000000000000000000000000000000000000000000000145884204447581481660582048048153974971515945684357636342943629684526  \\
				64    & -0.00000000000000000000000000000000000000000000000000000000000017667961728170528000576792667621868411541660759849978590782632609593 \\
				65    & 0.00000000000000000000000000000000000000000000000000000000000002140811101251326040856439038824440240421444446232730304347414207885  \\
				66    & -0.00000000000000000000000000000000000000000000000000000000000000259523933771640839395506227308797276439265244435789124567628998166 \\
				67    & 0.00000000000000000000000000000000000000000000000000000000000000031475844970368527207859988440756946780234602898258301523330011604  \\
				68    & -0.00000000000000000000000000000000000000000000000000000000000000003819199275948086232250066986474409024788436219755205198795932355 \\
				69    & 0.00000000000000000000000000000000000000000000000000000000000000000463613921141885306633052587623577878432147020683041743030130178  \\
				70    & -0.00000000000000000000000000000000000000000000000000000000000000000056302071895358497326551329253255067928964476053463704432078804 \\
				71    & 0.00000000000000000000000000000000000000000000000000000000000000000006840232492396395663308128089976226674456407391969807723301249  \\
				72    & -0.00000000000000000000000000000000000000000000000000000000000000000000831361414863747302543682188983858540128315562437846855903945 \\
				73    & 0.00000000000000000000000000000000000000000000000000000000000000000000101087762475684690072626638918750137063439086931801825227150  \\
				74    & -0.00000000000000000000000000000000000000000000000000000000000000000000012284024804151465169365011337828818556201671466330776051102 \\
				75    & 0.00000000000000000000000000000000000000000000000000000000000000000000001523094504376195525840058198589089909099106106553807071320  \\
				76    & -0.00000000000000000000000000000000000000000000000000000000000000000000000117299039745586326404875273916993576491873683755772837791 \\
				77    & 0.00000000000000000000000000000000000000000000000000000000000000000000000174592008510384125277934069832205616241191670358249124656  \\
				78    & 0.00000000000000000000000000000000000000000000000000000000000000000000000349419675898782669372040249072202278479470709387245197289  \\
				79    & 0.00000000000000000000000000000000000000000000000000000000000000000000000799670625368605793111544759252532945626113982364236937884  \\
				80    & 0.00000000000000000000000000000000000000000000000000000000000000000000001783736086252667800700690586334884589640117941023274188714  \\
				81    & 0.00000000000000000000000000000000000000000000000000000000000000000000003914629772923365981969623359065939745500344906392471974273  \\
				82    & 0.00000000000000000000000000000000000000000000000000000000000000000000008451862719318529215558189995798248228386015610623753894175  \\
				83    & 0.00000000000000000000000000000000000000000000000000000000000000000000017959488730997370412819250470216213014764808997616858075849  \\
				84    & 0.00000000000000000000000000000000000000000000000000000000000000000000037573129823027411320438637897746148414417808350459988510191  \\
				85    & 0.00000000000000000000000000000000000000000000000000000000000000000000077421493364609642333947794718727741247348087249217298429214  \\
				86    & 0.00000000000000000000000000000000000000000000000000000000000000000000157181555687935344929847398330810280467371631725554948612501  \\
				87    & 0.00000000000000000000000000000000000000000000000000000000000000000000314520476093280876644811508736958575816347686854506891900222  \\
				88    & 0.00000000000000000000000000000000000000000000000000000000000000000000620513768424215799858556499474177582182607773184311918853263  \\
				89    & 0.00000000000000000000000000000000000000000000000000000000000000000001207406287495543277451628598395555346789389222171356692562722  \\
				90    & 0.00000000000000000000000000000000000000000000000000000000000000000002317906407517644807430502076115241195672454732221700657294525  \\
				91    & 0.00000000000000000000000000000000000000000000000000000000000000000004391535644933142572262434185505603203399456343798168972287851  \\
				92    & 0.00000000000000000000000000000000000000000000000000000000000000000008213898261760714205027379987554661751530583907334785461102139  \\
				93    & 0.00000000000000000000000000000000000000000000000000000000000000000015171403728662894570055491378016484782263596816769234952817404  \\
				94    & 0.00000000000000000000000000000000000000000000000000000000000000000027680490412828741154710809310842123167257343965385307676036850  \\
				95    & 0.00000000000000000000000000000000000000000000000000000000000000000049902044928126690148948862379411085268486828266377908190620281  \\
				96    & 0.00000000000000000000000000000000000000000000000000000000000000000088916314010515705566014683642516706025000215138127727366871509  \\
				97    & 0.00000000000000000000000000000000000000000000000000000000000000000156632610036345267314420852986615817614367561071128915960553046  \\
				98    & 0.00000000000000000000000000000000000000000000000000000000000000000272858050807521788889043076064881355400620297213673392963699625  \\
				99    & 0.00000000000000000000000000000000000000000000000000000000000000000470174044807231982937532481326520156804492073444884756225621778  \\
				100   & 0.00000000000000000000000000000000000000000000000000000000000000000801601513689780222612781336297434728778902276505467370784946042  \\
			\end{longtable}
		}
\end{landscape}
\section{Python code to compute upper and lower bounds} \label{sec code}
Given precomputed values for \(\vec{a}\) and \(\vec{b}\) it is relatively simple to compute
(non-rigorous) upper and lower bounds for \(\nu_2^2\) using the \texttt{python} code and the \texttt{mpmath}
package. Here is an example that uses 16 coefficients for both \(\vec{a},\vec{b}\) and Richardson extrapolation~\cite{sidi2003} to
estimate the required infinite series. When run, the above script gives bounds on \(\nu_2^2\) precise to the first 32 digits.

\begin{displaymath}
	\triangleleft \diamond \triangleright
\end{displaymath}

{\scriptsize
\begin{lstlisting}[language=Python, numbers=left]
from mpmath import mp, besselj, fabs, factorial, inf, mpf, nstr, nsum, pi, power

mp.dps = 128  # digits of precision for the computation
epsilon = mpf("1e-50")  # tolerance for series extrapolation

half = mpf("0.5")
pi_on_two = pi() / mpf(2)
four_on_pi = mpf(4) / pi()
four_on_three = mpf(4) / mpf(3)


a_raw = [
    "0.9860206557159394456200552925916864935048606984660867659294396429",
    "0.01489770489081501339564402620613813180180708471357153405253588246",
    "-0.0010020514702611618623176011272831244046954374945414171879461225",
    "0.00009204407950297445040292960936251992007718734943876628562920974158",
    "-0.000009222749045097728828924770054951562677689395741936839825232504163",
    "0.0000009624829377769523619143359767714486583855960157862298060806209062",
    "-0.0000001030825466711237741281471720057151615014390669675059116945030138",
    "0.00000001125442999323010937747021511295873118119894549766096679782964655",
    "-0.000000001247526888343701534748343542688166570551788621599324821760315743",
    "0.0000000001399983570650130576391202478311227452665491901794575795215352259",
    "-0.00000000001587108251776074729790745823335658801904906951652187759074734197",
    "0.000000000001814579367249737396442135781162248916865862480501175051507144029",
    "-0.0000000000002089593636121966668255461641156413600243015932298164695606951635",
    "0.00000000000002421093838520663536817465987901271450247981787657530210660769827",
    "-0.000000000000002820080540615667295860418397000186709552077272721342675500717904",
    "0.0000000000000003300013142065405729012802718842249459836838193432161048014012172",
]
b_raw = [
    "1.02206573731609906209062285721046071122061035703720832517287196",
    "-0.02316346059073084779982636448754076645713917827003597676237278346",
    "0.001170550097383719246348345976716843011939144437093896137463616238",
    "-0.00007845159525437915088852263273588842719505793132135524091157766925",
    "0.000006101782370446285371463634245302883258983942949246427102319026463",
    "-0.0000005200974936973445795799002679008593946736768461635427163052152714",
    "0.00000004716048719323026047229288414394196590695665909728438524673409874",
    "-0.000000004471264809260967502221347731032458818596936687059973281105908718",
    "0.0000000004384319408026955946402611166292815653609028048460793198243935738",
    "-0.00000000004413778800416698015815314844080164423842416912926126118301283395",
    "0.000000000004538504357723005058673881111666831256023609333041237569856360479",
    "-0.0000000000004748627897491203459135138921486404308021033117402122537914927201",
    "0.00000000000005041269366573337296639115219231576962662698765250656699572447564",
    "-0.000000000000005418416187945572348485577704094735652600218337097929409831514064",
    "0.0000000000000005885883290033309163608068290402251544028088755111511925945235308",
    "-0.00000000000000006452865046019881731607090203710559895373215919872457713898102322",
]

# cast the coefficients to mpf's and then normalise them
tmp = [mpf(x) for x in a_raw]
s = sum(tmp)
a = [x / s for x in tmp]

tmp = [mpf(x) for x in b_raw]
s = sum(tmp)
b = [x / s for x in tmp]  # we are using sum b = 1.


# We multiply the BesselJ functions by these amplitudes
j_amps = [factorial(p) * power(four_on_pi, p) for p in range(len(a) + 2)]


def big_J(p, k):
    return j_amps[p] * besselj(p, k * pi_on_two) / power(k, p)


def big_F(k):
    return power(sum([big_J(p, k) * v for p, v in enumerate(a)]), 4)


def big_G(k):
    return power(
        fabs(sum([big_J(p + 1, k) * v for p, v in enumerate(b)])), four_on_three
    )


# group the summands in 4s like this so that their asymptotics are better behaved.
# gives better series extrapolations
def the_F_summand(n):
    return big_F(4 * n - 3) + big_F(4 * n - 2) + big_F(4 * n - 1) + big_F(4 * n)


def the_G_summand(n):
    return big_G(4 * n - 3) + big_G(4 * n - 2) + big_G(4 * n - 1) + big_G(4 * n)


# use richardson extrapolations to estimate the infinite series
sum_F = nsum(the_F_summand, [1, inf], method="r", tol=epsilon, verbose=True)
sum_G = 2 * nsum(the_G_summand, [1, inf], method="r", tol=epsilon, verbose=True)

print("Sum of F(k)^4:\n\t", nstr(sum_F, 64))
print("Sum of abs(G(k))^(4/3):\n\t", nstr(sum_G, 64))

nu_u = half + sum_F
nu_l = half + half / power(sum_G, 3)

print("This gives bounds:")
print(nstr(nu_l, 36), "<= nu_2^2 <=")
print(nstr(nu_u, 36))
print("delta:", nstr(nu_u - nu_l, 4))
\end{lstlisting}
}
\begin{displaymath}
	\triangleleft \diamond \triangleright
\end{displaymath}

{\footnotesize
\begin{lstlisting}[]
----------------------------------------------------------------------
Adding terms #0-#10
Direct error: 0.000105053
Richardson error: 0.0746387
----------------------------------------------------------------------
Adding terms #10-#30
Direct error: 1.10602e-5
Richardson error: 9.5303e-7
----------------------------------------------------------------------
Adding terms #30-#60
Direct error: 2.72974e-6
Richardson error: 4.77855e-22
----------------------------------------------------------------------
Adding terms #60-#100
Direct error: 9.77732e-7
Richardson error: 7.45014e-51
----------------------------------------------------------------------
Adding terms #0-#10
Direct error: 0.00132661
Richardson error: 0.942536
----------------------------------------------------------------------
Adding terms #10-#30
Direct error: 0.000139668
Richardson error: 1.20349e-5
----------------------------------------------------------------------
Adding terms #30-#60
Direct error: 3.44711e-5
Richardson error: 6.03435e-21
----------------------------------------------------------------------
Adding terms #60-#100
Direct error: 1.23468e-5
Richardson error: 9.39307e-50
----------------------------------------------------------------------
Adding terms #100-#150
Direct error: 5.47363e-6
Richardson error: 2.02238e-85
Sum of F(k)^4:
	 0.07463960715151959272725542752705485435561425481480583579118891616
Sum of abs(G(k))^(4/3):
	 1.88509635262914968733569626740377511940895683024823141154615938
This gives bounds:
0.574639607151519592727255427527052952 <= nu_2^2 <=
0.574639607151519592727255427527054854
delta: 1.903e-33
\end{lstlisting}
}

\begin{displaymath}
	\triangleleft \diamond \triangleright
\end{displaymath}

\end{document}